\documentclass{amsart}
\usepackage{graphics}
\usepackage{verbatim}
\def\scorept{\text{ \sl pt}}
\def\qrtet{{q.r.~tet}}
\def\qrtets{{q.r.~tets}}

\theoremstyle{plain}% default
\newtheorem{lem}{Lemma}[section]
\newtheorem{thm}{Theorem}[section]

\theoremstyle{remark}
\newtheorem{calc}{Calculation}[subsection]
\newtheorem*{remark}{Remark}

\DeclareMathOperator{\score}{sc}
\DeclareMathOperator{\sol}{sol}
\DeclareMathOperator{\dih}{dih}
\DeclareMathOperator{\vor}{vor}
\DeclareMathOperator{\octavor}{octavor}
\DeclareMathOperator{\gma}{gma}

\begin{document}

\begin{abstract}

The Hales program to prove the Kepler conjecture on sphere packings
consists of five steps, which if completed, 
will jointly comprise a proof of the conjecture.
We carry out step five of the program, a proof that the local density of a
certain combinatorial arrangement, the pentahedral prism, is less
than that of the face-centered cubic lattice packing.
We prove various relations on the local density using computer-based
interval arithmetic methods.  Together, these relations imply the
local density bound.

\end{abstract}

\title[Sphere Packings, V]{Sphere Packings, V}
\author{Samuel L. P. Ferguson}
\address{Department of Mathematics, University of Michigan, Ann Arbor,
MI 48109}
\email{samf@math.lsa.umich.edu}
\urladdr{http://www.math.lsa.umich.edu/\~{}samf}
\date{\today}
\maketitle

\section{Introduction}
A collection of uniformly-sized balls in Euclidean 3-space is called
a {\em packing} if no two balls have a common interior point.  We
refer to such a packing as a {\em sphere packing}.
The Kepler conjecture asserts that the density of a packing of equal
spheres in three-dimensions cannot exceed that of the face-centered
cubic packing.
Recently, Thomas C. Hales proposed a program designed to prove the
Kepler conjecture \cite{remarks, problem, sp1}.  
This program proposes a formulation designed
to reduce the complexity of the problem to a level tractable via 
modern methods.

This formulation requires first that we define a decomposition of
space, relative to a given packing of spheres.
Next, we define stars, based on the decomposition.
We then define a local density for these stars.
Finally, we show that a sufficient local density bound for stars implies
the required global density bound for the packing.

In addition to the formulation of the conjecture,
the Hales program consists of five steps.
An overview of the program detailing the content of each step
can be found in \cite{sp1}.
Details of the formulation itself can be found in \cite{FKC}.
We carry out step five of the program.

%Decomposition of space
\subsection{A Decomposition of Space}
A decomposition of space, relative to a given packing of spheres,
is critical to the formulation of the program.
Our decomposition is inspired by the Delaunay and Voronoi decompositions,
drawing from the advantages of each.
We define a decomposition
of space composed of a {\em $Q$-system} together with modified Voronoi
cells called {\em $V$-cells}.  
We first define the $Q$-system, a partial decomposition of space,
and then complete the decomposition by adding the $V$-cells.

\begin{remark}
There are many possible decompositions which could be used in a
similar approach to the conjecture.  The challenge is to select
a decomposition which brings each element of the proof within reach.  
This selection has not been at all easy.  For example, both
the Delaunay and Voronoi decompositions have critical defects which
prevent their use.  The $Q$-system is so named because it consists
of certain tetrahedrons called \qrtets\ and quarters.
\end{remark}

We say that two sets {\em overlap} if their relative interiors
intersect.
To define the $Q$-system, we begin with an arbitrary saturated packing 
of spheres of radius one.  
We say that a packing is {\em saturated} if no more spheres may be
added to the packing without overlapping spheres already in the packing.  
We consider only saturated packings, 
as we are looking for packings of maximal density.
We identify the centers of the spheres as
{\em vertices} of the packing.  On occasion, we will distinguish
a vertex of the packing, calling it the {\em origin}.

\subsubsection{Identification by edge lengths}
\begin{figure}
\includegraphics{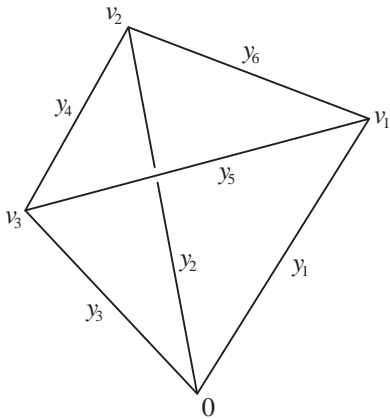}
%\centerline{ \psfig{file=tet.eps} }
\caption{Tetrahedron with distinguished vertex and labeled edges.}
\label{fig:tet}
\end{figure}

By connecting vertices in the
packing, we can identify certain tetrahedrons.  
This identification is based on the
length of the edges connecting the vertices.
We label the edges and vertices of tetrahedrons 
as in Figure~\ref{fig:tet}.  Note that the order of
the edge lengths determines which is the distinguished vertex
of a tetrahedron.

\subsubsection{Quarters and \qrtets}
We define two types of tetrahedrons, a 
{\em quasi-regular tetrahedron}, called a {\em \qrtet}, and a
{\em quarter}.
A \qrtet\ is a tetrahedron whose edge lengths each lie
in the interval $[2,2.51]$.  A quarter is a tetrahedron with five
edges with lengths in $[2,2.51]$ and one long edge
 with length in $[2.51,2 \sqrt{2}]$.
 
If a vertex of
a quarter is distinguished, the quarter has
two orientations:  flat, and upright.  A quarter is said
to be {\em flat} if the distinguished vertex is opposite the long
edge.  Similarly, a quarter is said to be {\em upright} if
the distinguished vertex lies along the long edge.  We refer to
the long edge of a quarter as the {\em diagonal}.

\begin{remark}
Throughout this paper we introduce various special constants, such
as $2.51$.  Although they may appear to be rather mundane numbers,
these constants have been carefully selected.  It is usually the
case that our selection comes from a range of numbers.  We tend
to choose numbers which are easily represented in decimal form,
although for some purposes, numbers with a simple binary representation
would be better.  We do not mean to imply that our computations
are made to only two-digit significance, nor that the selection
of constants is arbitrary.  
\end{remark}

%The $Q$-system is composed of a non-overlapping collection of
%\qrtets\ and quarters.  The selection
%of \qrtets\ and quarters is somewhat complicated, however.

\subsubsection{Octahedra}
If four quarters fit together along a common diagonal, forming
a figure with six vertices, the resulting figure 
composed of quarters is called
an {\em octahedron}.  We use the term octahedron only when we
refer to such a figure composed of quarters.
An octahedron may have more than one diagonal
of length at most $2\sqrt{2}$, so the decomposition of an
octahedron into quarters may not be unique.

%Define a quad cluster?  Or do it later?

\subsubsection{The $Q$-system}
The $Q$-system is a collection of non-overlapping \qrtets\ and
quarters.
We construct the $Q$-system incrementally, building it from
\qrtets\ and quarters.  To begin, we identify all \qrtets,
quarters, and octahedra in the packing.  We will define
the $Q$-system to be a subset of these structures.
We consider each \qrtet\ and quarter in turn, choosing which
ones we will add to the $Q$-system.

We first consider all identified octahedra.
For each octahedron in the packing, we fix a diagonal of length at most 
$2\sqrt{2}$, and place the four quarters along that diagonal
in the $Q$-system.  
Next, we place all quasi-regular tetrahedra
in the $Q$-system.  By Lemmas~1.2 and 1.3 of \cite{FKC}, 
it is not possible for \qrtets\ to
overlap either themselves or octahedra.

Next, we wish to place the rest of the quarters into
the $Q$-system.  Unfortunately, there is some ambiguity about how to do this
without overlapping the quarters.

\begin{remark}
The resolution of the ambiguities in placing quarters into the
$Q$-system is tedious.  Details can be found in \cite{FKC}.
\end{remark}

We say that two tetrahedra are {\em adjacent} if
they share a face.

The first ambiguity arises when we have two adjacent
quarters, sharing a diagonal.  If we identify their common vertex
opposite the diagonal, they become flat quarters.  It is possible that
the same collection of vertices could be decomposed into two flat
quarters with a different diagonal.  We call this ambiguity a
conflicting diagonal.  It is not important for
this paper which choice we make for the decomposition.  Making a
selection when necessary, we place the adjacent quarters which
share a diagonal into the $Q$-system.

The second ambiguity which arises is the case of two isolated,
overlapping quarters.  By isolated we mean that neither is
adjacent to another quarter with which it shares a diagonal.
If such a case arises, we place neither quarter in the $Q$-system.

To complete the $Q$-system, we add all remaining unconsidered quarters, 
those which are isolated
and which do not overlap anything else in the $Q$-system.
This completes the definition of the $Q$-system.

\subsubsection{Voronoi cells}
The {\em Voronoi cell} at a vertex of the packing consists of
all points in space which are closer to that vertex
than any other vertex of the packing.  The Voronoi
decomposition of space is simply the collection of all
Voronoi cells.  

\begin{remark}
To define a $V$-cell at a vertex $w$,
we could take the intersection of the
Voronoi cell at $w$ with the complement of the $Q$-system.
However, this formulation
leads to unnecessarily complex $V$-cells.  
This is due to the fact that a tetrahedron
need not contain its circumcenter.
\end{remark}

In order to define $V$-cells, we must first discuss the
orientation of a vertex with respect to a face.

If the circumcenter of a tetrahedron $S$ lies on the side 
of the plane $P$ through the face of $S$ opposite the vertex,
we say that the face has {\em negative orientation} with
respect to that vertex.  To identify points lying on the
side of $P$ opposite the vertex, we say that these points
lie on the {\em negative side} of the face.  
By Lemma~2.1 of \cite{FKC}, at
most one face of a quarter or \qrtet\ has negative
orientation.  

\begin{remark}
If a face of an element $S$ of the $Q$-system
has negative orientation, the Voronoi cell associated
with the vertex of $S$ opposite that face will pass through the
face, unless the Voronoi cell is truncated by a
vertex outside of $S$.
Hence the tip of the Voronoi
cell protruding through the tetrahedron would be ``orphaned"
in the sense that the tip would be part of the star
but not be contiguous with the rest of the Voronoi cell, since
it will be separated by an element of the $Q$-system.
By design, $V$-cells reapportion the protruding
tip into the adjacent Voronoi cells so that we may 
avoid this complication.
\end{remark}

\subsubsection{$V$-cells}
The construction of the $V$-cells is somewhat complex.
We first describe a vertex deletion process, which we
will use to construct the $V$-cells.  

For each negatively-oriented
tetrahedron $S$, we distinguish the vertex $v$ opposite the
negatively-oriented face.
%In the deletion process,
%we view each negatively-oriented tetrahedron $S$
%in the $Q$-system in isolation, meaning that we ignore the
%effect of all other vertices in the packing on the Voronoi
%cell determined by $S$.

\begin{figure}
\includegraphics{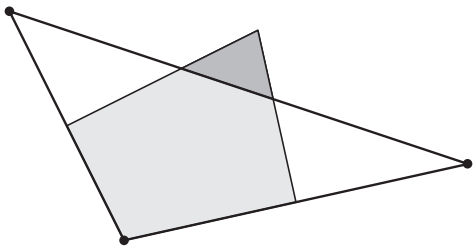}
%\centerline{ \psfig{file=twodimvor.eps} }
\caption{Two dimensional simplex with Voronoi cell.}
\label{fig:twodimvor}
\end{figure}

We define the {\em tip} of the Voronoi cell associated with
an isolated negatively-oriented tetrahedron $S$
to be the part of the Voronoi cell determined by $S$ alone
which lies on the negative side of the face opposite $v$.
Figure~\ref{fig:twodimvor} represents the two-dimensional
analog of a tip.

%If we view each negatively-oriented tetrahedron
%in isolation, each point
Each point
$x$ outside the $Q$-system lies in finitely many tips,
each associated with a simplex $S$ with distinguished
vertex $v$.  (Typically, this collection of vertices will be empty.)
Temporarily deleting these vertices from the packing, we take the
Voronoi decomposition of the remaining collection of
vertices.  The point $x$ lies in the modified Voronoi
cell at some new vertex $u$.  

The deletion process therefore
provides a mapping from the complement of the $Q$-system to
the vertices of the packing.

We now consider a vertex $w$ of the packing.  We define
the $V$-cell at $w$
to be the set of points $x$ outside the $Q$-system
with the property that each lies in a
modified Voronoi cell at $w$ after the deletion process.  That
is, the $V$-cell at $w$ is the pre-image of $w$ under the
mapping.

The collection of \qrtets\ and quarters in the $Q$-system together with
all the $V$-cells comprises our decomposition of space, 
relative to a fixed packing of spheres.

%Define standard regions.
\subsection{Decomposition Stars}
Now that we have defined a decomposition of space relative to
a packing, we can define a decomposition star relative to a
vertex of the packing.  We introduce the concept of a cluster
of vertices, and then define a decomposition star in terms of
the graph of the cluster.

We define the {\em height} of a vertex to be its distance
from the origin.

Fix the origin at a vertex of the packing.
Mark all vertices whose height is at most
$2.51$.  The collection of marked vertices is called a {\em cluster}.

If the distance between two vertices in the cluster does not
exceed $2.51$,
join them by an edge.  We call the resulting structure the
{\em graph} of the cluster.

Next, project the graph radially to the unit sphere centered at the
origin.
This projection produces a partition of the surface
into regions called {\em standard regions}.
By Lemma~3.10 of \cite{sp1}, each edge projects to an
arc on the unit sphere, and these arcs do not meet except at
endpoints.

We say that a vertex is {\em enclosed} over
a figure if the vertex lies in the infinite cone at the origin 
generated by the figure.

We associate a {\em decomposition star} with each vertex $v$.
We define a decomposition star
as the union of the $V$-cell at $v$ with all the
simplices in the $Q$-system having a vertex at $v$.  

Decomposition stars and the standard regions have several
important properties.

In the most general terms, we wish to show that
$V$-cells, the $Q$-system and decomposition stars are compatible
with standard regions.  In the discussion which follows, we
clarify this statement.

We will define two
properties, each of which depends on the local
nature of $V$-cells and the $Q$-system.  These two
properties define the compatibility of $V$-cells and the
$Q$-system with standard regions.

The first property is that each
\qrtet\ or quarter of a decomposition star 
lies inside the cone over a single standard region.

The second property is somewhat harder to define.

We define a {\em quasi-regular face} to be a triangle 
constructed from vertices in the packing, whose edge
lengths lie in the interval $[2,2.51]$.

Note that a side of a standard region corresponds 
to a quasi-regular face.

By Lemma~2.2 of \cite{FKC}, if a vertex has negative
orientation with respect to a quasi-regular face, it must form
a \qrtet\ with that face.

Using this information, we now consider the 
intersection of a $V$-cell with the cone over a
standard region.  We claim that the shape of the
$V$-cell inside the cone is completely determined 
by the vertices lying inside the cone.  

If a vertex $w$ outside the cone 
is to affect the $V$-cell inside the cone,
it must have negative orientation with respect to
a quasi-regular face $F$ corresponding to a side of the
standard region.  By the previous result, the vertex $w$
must form a \qrtet\ with a face of the standard region.
See 2.2 of \cite{sp2}.

By construction, if the Voronoi cell at $w$ protrudes
through $F$, the protruding part is attached to the $V$-cells
inside the cone.
%fix this.  
Therefore the $V$-cell lying inside the
cone over a standard region is completely determined by the
vertices lying inside the cone.

\begin{remark}
These properties of standard regions allow us to consider each
standard region in isolation.  If this were not possible, the
treatment of each decomposition star would be hopelessly complex.
\end{remark}

\subsubsection{Standard clusters}
In our investigation of decomposition stars, we will need to
consider the geometric arrangements associated with all possible
standard regions.  We call these arrangements standard clusters.
We consider standard clusters in the context of decomposition stars.

Specifically, the {\em standard cluster} associated 
with a standard region $R$ of a decomposition star
is the union of the simplices in the decomposition star
which lie in the cone over $R$ together with the
part of the $V$-cell that lies over $R$.

The standard cluster associated with a triangular standard
region is a \qrtet.

We call the standard cluster associated with a
quadrilateral standard region a {\em quad cluster}.
%For example, an octahedron is a quad cluster, as are two
%adjacent flat quarters joined along their diagonal.

%Define decomposition stars and local structures.
%Quoting, more or less, from FKC . . .

%Classification of quad clusters.
\subsubsection{Classification of Quad Clusters}
\label{quad:classification}
Our treatment of decomposition stars requires that we
classify all possible quad clusters.

\begin{remark}
We single out quad clusters because the decomposition star
which we treat in this paper has only triangular and
quadrilateral standard regions.
\end{remark}

We call the four vertices of the packing
which project to the vertices of the
quadrilateral region the {\em corners} of the quad cluster.

We give two exhaustive lists of the possible decompositions,
based on the length of the diagonals between the corners
of a quad cluster.

If a quad cluster has a diagonal between two corners
whose length does not exceed
$2\sqrt{2}$, there are three possible decompositions.  

First,
the quad cluster could be composed of two flat quarters sharing
the diagonal.  We call this a {\em flat} quad cluster.
Second, the quad cluster could be composed of
two flat quarters sharing the other diagonal.  Third, the
quad cluster could be composed of four upright
quarters forming an octahedron.

If both diagonals between the corners of a quad cluster
have lengths which exceed $2\sqrt{2}$, there are again
three possible decompositions.

First, the quad cluster could be an octahedron.  Second,
there could be no enclosed vertex of height at most $2\sqrt{2}$.
We call such a quad cluster a {\em pure Voronoi} quad cluster,
as it consists of only the part of the 
$V$-cell lying above the quadrilateral,
since such a quad cluster cannot contain any quarters.
Third, there could be an enclosed vertex of height at most
$2\sqrt{2}$.  We call such arrangements {\em mixed} quad
clusters, since they consist of zero or more upright quarters and 
portions of $V$-cells lying above the quadrilateral.

\subsection{A Local Density Function}
The next element in our formulation is a local density function
for decomposition stars.

We first construct a local density function for standard clusters.
We will construct a local density function for decomposition
stars by applying the local density function to each standard
cluster in the star.

%We will define a local density function on
%decomposition stars.  An upper bound on the local density function
%for each decomposition star will imply an upper bound on the global
%density of a packing.
%Add the argument here?  Or wait until later?

%Put this in the section on bounds?
\begin{comment}
Due to the complex nature of the local density function, it is typically
unrealistic to attempt to prove the required bounds directly.  Instead,
we prove the
majority of the required bounds via
computer-based interval arithmetic methods,
while taking advantage of various reduction arguments to reduce
the complexity of the verifications.
\end{comment}

\subsubsection{The score of a standard cluster}
The {\em score} of a standard cluster represents
one of many possible local density functions applied to a
standard cluster.  The choice of function 
depends on the geometry of the standard cluster in question. 
Our selection is made in an attempt to produce the best possible bounds
on the local density.

In this paper we construct scoring functions only for
\qrtets\ and quad clusters.  The scoring functions for other
standard clusters can be found in other papers in the program.

We further subdivide
our scoring functions, applying them to component quarters of
quad clusters.  

We define functions on tetrahedrons in terms
of the edge lengths $(y_1,\ldots,y_6)$ of a tetrahedron.  The
first three edges $y_1,y_2,y_3$ are adjacent to the
distinguished vertex or origin.  Edge $y_i$ lies opposite
edge $y_{i+3}$ for $i=1,2,3$.

We use several distinct scoring functions.  These functions are
typically linear combinations of Voronoi volume and solid angle.
To define these scoring functions, we require some intermediate
definitions.

We define the solid angle {\em sol(S)} of a
tetrahedron $S$ with distinguished vertex $v$ 
to be three times the volume of the intersection
of the unit ball at $v$ with the infinite cone on $S$
with origin at $v$.  The units of solid angle are
therefore steradians.

We define the geometric Voronoi volume of a tetrahedron $S$ with distinguished
vertex $v$ to be the volume occupied by points of $S$ lying closer
to $v$ than to the other vertices of the tetrahedron.

We define the analytic Voronoi volume $\text{vol}(S)$ 
of a tetrahedron $S$ with distinguished vertex $v$ to be
the analytic continuation of the formula for the geometric Voronoi
volume which holds when $S$
contains its circumcenter.  The analytic Voronoi volume is
identical to the geometric Voronoi volume only when $S$ contains
its circumcenter.  However, the sum of the analytic Voronoi
volumes at all four vertices of $S$ is equal to the volume of $S$.

\begin{remark}
This property allows us to produce a global density bound from a bound
on the score of a decomposition star.
\end{remark}

 Finally, let
$\delta_{\text{oct}} = (\pi - 4\arctan(\sqrt{2}/5))/(2\sqrt{2})$.

\subsubsection{The score of a tetrahedron}
We define three scoring functions for tetrahedrons.

The simplest is called
vor analytic or {\em vor}, taking its name from the analytic Voronoi volume.
The vor analytic scoring $\vor(S)$ of a tetrahedron $S$ is given by
$4(-\delta_{\text{oct}} \text{vol}(S) + \sol(S)/3)$.  Here {\em vol}
denotes the analytic Voronoi volume at the distinguished vertex $v$.

Averaging $\vor(S)$ over the two vertices adjacent to the long edge
of an upright quarter gives {\em octavor}.  That is, if the long edge
of an upright quarter is $y_1$, 
\[
\octavor(S) = (\vor(y_1,y_2,y_3,y_4,y_5,y_6) + \vor(y_1,y_5,y_6,y_4,y_2,y_3))/2.
\]

The average of $\vor(S)$ over all four vertices of a
tetrahedron gives {\em gma}, which we also call the {\em compression} of
a tetrahedron.  That is,

\begin{align*}
\gma(S) = & (\vor(y_1,y_2,y_3,y_4,y_5,y_6) + \vor(y_1,y_5,y_6,y_4,y_2,y_3)  \\
	& + \vor(y_2,y_4,y_6,y_5,y_1,y_3) + \vor(y_3,y_4,y_5,y_6,y_1,y_2))/4.
\end{align*}

Our method for choosing the scoring function for a particular
tetrahedron depends on whether it is a \qrtet\ or a quarter.

\subsubsection{Scoring a \qrtet}
If the circumradius of a \qrtet\ $S$ is less than $1.41$, we score
$S$ by compression.  Otherwise, we score $S$ using
vor analytic.

\subsubsection{Scoring a quarter}
\label{scoring:quarter}
For quarters, the scoring system depends on both the
orientation of the quarter and the circumradii of the faces adjacent to
the diagonal.  If the circumradius of each of the adjacent
faces to the diagonal is less than $\sqrt{2}$, the quarter is scored using
compression.  Otherwise, flat quarters are scored using
vor analytic, and upright quarters are scored using 
{\em octavor}.

\subsubsection{The score of $V$-cells}
The last element to be scored is the 
intersection of a $V$-cell with the cone over a
standard region.  Call this intersection $S$.
Such elements are scored using
Voronoi scoring.  That is, they are scored using the usual formula
$4(-\delta_{\text{oct}} \text{vol}(S) + \sol(S)/3)$, 
where in this case $\text{vol}(S)$ represents the volume of $S$,
and $\sol(S)$ represents three times the volume of the intersection of
$S$ with the unit ball at the origin.

\subsubsection{Truncated Voronoi scoring}
If a quad cluster contains only vertices whose
height exceeds $2 \sqrt{2}$, we may bound the score of the
quad cluster by using
{\em truncated Voronoi}, meaning
that the $V$-cell is truncated at a distance of
$\sqrt{2}$ from the distinguished vertex.  As this
method decreases the volume term, it provides an upper bound
on the Voronoi score.
The bound on scoring
is then $4(-\delta_{\text{oct}}\text{vol} + \sol/3)$, where
in this case, $\text{vol}$ denotes the truncated Voronoi 
volume associated with the distinguished vertex, and
$\sol$ denotes the solid angle associated with the
distinguished vertex.

\begin{remark}
The selection of the scoring functions was a challenging part of
the formulation of the problem.  Compression
arose as a natural local density function for the Delaunay decomposition.
Similarly, $\vor$ arose as a local density function for the
Voronoi decomposition.  The scoring formulation which we exhibit
is the result of much experimentation.
\end{remark}

\subsubsection{Summary of functions}
We use the following functions:  {\em sol},
the spherical angle associated
with the distinguished vertex of a tetrahedron or quad cluster, 
{\em dih}, the
dihedral angle associated with the first edge of a tetrahedron,
{\em vor}, the vor analytic score of a tetrahedron,
{\em octavor}, the averaged vor analytic score of an
upright quarter,
{\em gma}, the compression of a tetrahedron,
and {\em sc}, a generic name for
the score of a tetrahedron or quad cluster.

Explicit formulas for these special functions are available in
\cite{sp1, sp2, FKC}.

\begin{remark}
We realize that our naming conventions are somewhat cumbersome.
Initially, we gave traditional symbols to the special functions.
For example, we referred to $\gma$ as $\Gamma$.  
However, as the library of functions increased, it became 
simpler to give each function a name which could be written easily
in both {\sl Mathematica} and C.  Greek symbols do not
have this property.
\end{remark}

\subsubsection{The score of a decomposition star}
The score of a decomposition star is the sum of the scores
of the standard clusters which comprise the
star, such as \qrtets\ and quad clusters.

\subsubsection{Definition of a point}
We find it useful to define a {\em point} 
to be the score of a regular tetrahedron (whose edge lengths
are $(2,2,2,2,2,2)$),
%$\scorept = 4\arctan(\sqrt{2}/5) - \pi/3$.
\[
\scorept = 4\arctan(\frac{\sqrt{2}}{5}) - \frac{\pi}{3}.
\]

\subsection{The global density bound}
The final part of our formulation of the Kepler conjecture
requires that we demonstrate
how a scoring bound on all decomposition stars implies a global
bound on the density of a packing.

Our formulation expresses the required
global density bound of $\pi/\sqrt{18}$ in terms of a proposed bound on
the score of a decomposition star.  Together, the five steps in
the Hales program will imply that
the score of any decomposition star does
not exceed $8 \scorept$.

In order to produce the transition between local and global
density bounds, we first discuss the essential equivalence of
all of the scoring methods.

Recall that the analytic Voronoi decomposition of a tetrahedron
produces a partition with four pieces.  We compute the
analytic Voronoi score of a tetrahedron using the analytic
continuation of the volume of the
piece containing the distinguished vertex.  Adding these four
(analytically continued) volumes gives the volume of the tetrahedron.

If we average $\score$ over all vertices of a tetrahedron $S$, regardless of
the scoring method, we arrive at $\gma(S)$.  In this sense, the
different scoring schemes for tetrahedra are all equivalent.
Recall that $\gma(S)$ is a linear combination of the volume of $S$
and the sum of the solid angles at each vertex of $S$.  That is,
\[
\gma(S) = -\delta_{\text{oct}} \text{vol}(S) + \frac{1}{3}\sum_{i=1}^4 \sol_i,
\]
where $\sol_i$ represents the solid angle of $S$ at vertex $i$.

Recall that the Voronoi score of a $V$-cell at a vertex $w$ is
\[
4(-\delta_{\text{oct}} \text{vol} + \frac{1}{3}\sol),
\]
where $\text{vol}$ represents the volume of the $V$-cell and
$\sol$ represents the solid angle at $w$ of the $V$-cell.

Begin with a sphere of radius $R$, containing $N$ vertices.
Assume that the proposed bound,
\[
\score \le 8 \scorept,
\]
holds for every decomposition star.
Add this bound for every star in the sphere.  Call the sum of the
scores $\sum \score$.  Each tetrahedron is shared by four stars.
Therefore the contribution to $\sum \score$ from each tetrahedron $S$
is $4 \gma(S)$.

The volume terms add up to $\frac{4 \pi}{3} R^3$ (neglecting the
boundary), and the solid angles add to give $\frac{4 \pi}{3} N$,
the volume of the unit balls within the sphere of radius $R$.
The boundary term is negligible as $R\to\infty$.

This gives, up to a negligible boundary term,
\[
4(-\delta_{\text{oct}} \frac{4 \pi}{3} R^3 + \frac{4 \pi}{3} N) \le 8 N \scorept.
\]
The density of balls within the sphere (again, up to a negligible boundary
term) is then
\[
\frac{\frac{4 \pi}{3} N}{\frac{4 \pi}{3} R^3} \le 
	\frac{2 \pi \delta_{\text{oct}}}{2 \pi - 3\scorept}
\]

Simplifying the right hand side, we find that it is equal to the desired
bound, $\pi/\sqrt{18}$.  Therefore the density of any packing cannot
exceed the conjectured bound of $\pi/\sqrt{18}$.

For a more detailed account of this computation, see \cite{sp1,FKC}.

\begin{figure}
\includegraphics{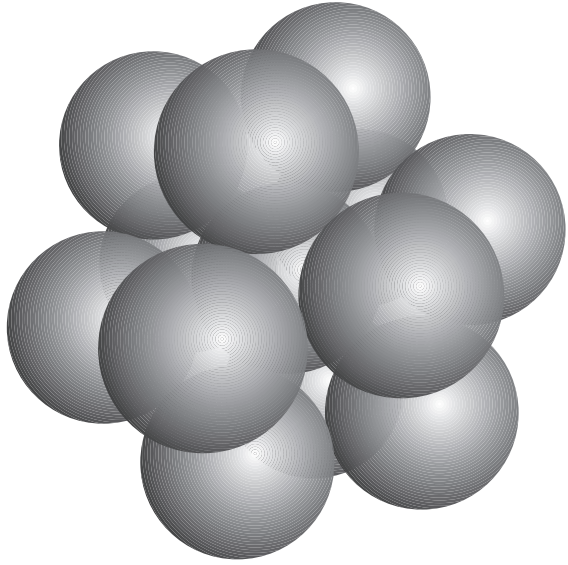}
%\centerline{ \psfig{file=fccballs.eps} }
\caption{The face-centered cubic cluster.}
\label{fig:fccballs}
\end{figure}

\begin{figure}
\includegraphics{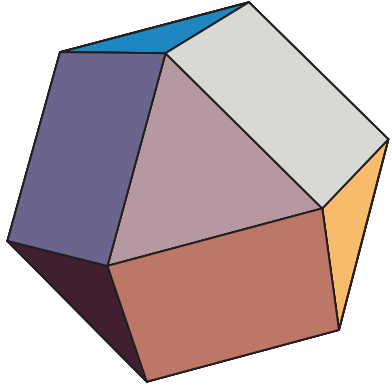}
%\centerline{ \psfig{file=fccface.eps} }
\caption{The face-centered cubic decomposition star.}
\label{fig:fccface}
\end{figure}

\begin{remark}
The proposed $8 \scorept$ bound is achieved by the decomposition star
of the face-centered cubic lattice packing.  See Figure~\ref{fig:fccballs}.
This star is composed of
8 regular \qrtets\ and 6 regular quad clusters.  See Figure~\ref{fig:fccface}.
Each \qrtet\ in this
star scores $1 \scorept$, the maximum possible score of a \qrtet.
Each quad cluster scores $0 \scorept$, the maximum score of a quad
cluster.  Hence the score of a face-centered cubic star is $8 \scorept$.
\end{remark}

%\subsection{The Hales program}
%Now that we have presented our formulation of the Kepler conjecture,
%we can explain the steps of the Hales program.

\section{The Pentahedral Prism}
Having introduced our formulation of the Kepler conjecture, we can
now define the subject of this paper, the pentahedral prism.

\begin{figure}
\includegraphics{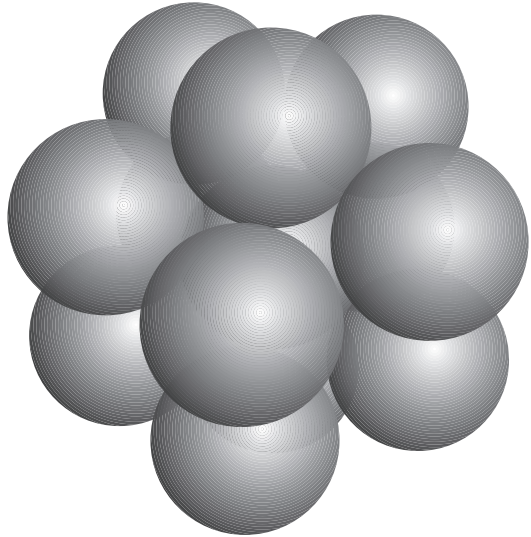}
%\centerline{ \psfig{file=pentaballs.eps} }
\caption{The pentahedral prism.}
\label{fig:penta2}
\end{figure}

The pentahedral prism arises as a decomposition star 
in a saturated packing of 3-space
with spheres of unit radius.
It is a particular cluster with twelve vertices (not counting the origin).
See Figure~\ref{fig:penta2}.

\begin{figure}
\includegraphics{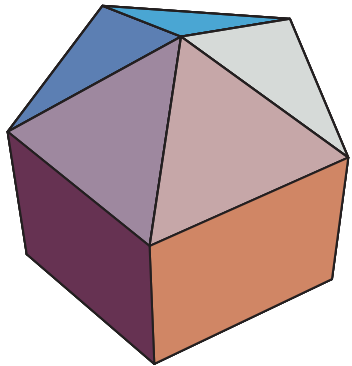}
%\centerline{ \psfig{file=pentaface.eps} }
\caption{The faces of a pentahedral prism.}
\label{fig:penta}
\end{figure}

The pentahedral prism is characterized by the arrangement and
combinatorics of its standard regions.  It is composed of ten triangular
standard regions, and five quadrilateral standard regions.
The ten triangles are arranged in two
{\em pentahedral caps}, five triangles arranged around a common vertex.
The five quadrilaterals lie in a band between the two caps.
See Figure~\ref{fig:penta}.

Recall that the standard cluster attached to a triangular standard region
is a \qrtet.  Likewise, the standard cluster attached to a 
quadrilateral is a quad cluster.  
We use the term pentahedral cap to refer to
both the standard regions and the \qrtets\ which comprise it.

The five steps in the Hales program \cite{sp1} are
\begin{enumerate}
\renewcommand{\theenumi}{\Roman{enumi}}
\item
A proof that even if all standard regions of a decomposition star are 
triangular the score is less than $8 \scorept$.
\item
A proof that standard clusters in regions of more than three sides
score at most $0 \scorept$.
\item
A proof that if all of the standard regions of a decomposition star are
triangles or quadrilaterals, then the score does not exceed $8 \scorept$
(excluding the case of the pentahedral prism).
\item
A proof that if some standard region has more than four sides, then
the star scores less than $8 \scorept$.
\item
A proof that the pentahedral prism scores less than $8 \scorept$.
\end{enumerate}

\begin{remark}
These steps were originally intended to be of roughly equal magnitude.
Due to changes in our formulation of the decomposition of space and
the local density function, some steps have become easier to complete
than others.  

In the earliest formulation, the pentahedral prism
could score more than $8 \scorept$.  It appeared, on an experimental
basis, to be the sole counterexample to the potential success of
the program.  To a certain degree, the pursuit of a proof of step
five drove the evolution of both the decomposition of space and the
local density function.  

The pentahedral prism still comes remarkably
close to achieving the optimal local density, that achieved by the
decomposition stars of the face-centered cubic lattice packing, which
score $8 \scorept$.
In this sense, we consider the pentahedral prism 
to be a ``worst case" decomposition
star.  As such, it required the devotion of much time and resources.
Techniques which we developed to handle this case have proven to
have significant impact on other steps of the program (as, conversely,
techniques from other steps have affected this one).  

Thus far,
the relations required to treat the pentahedral prism have been
delicate in contrast to the more general bounds which have to date
sufficed to treat other decomposition stars.
\end{remark}

As of this writing, steps I and II are complete \cite{sp1, sp2}.
Hales has exhibited partial results for III and IV.  

\begin{remark}
The nomenclature for the star which we treat in step five is
somewhat problematic.  This object is a prism in only the weakest 
sense, technically
speaking.  In addition, it is denoted a ``pentagonal prism"
in \cite{sp1}.  With a similar degree of inaccuracy, we could also
call it a ``pentangular prism".
\end{remark}

\section{The Scoring Bound}
\label{sec:dabound}
We present computations using auxiliary bounds 
which imply the main result of the paper, that the
score of the pentahedral prism is less than $8 \scorept$.

We use the following lemma in proving almost all of the
auxiliary bounds.

\begin{lem}
\label{lem:pentscorebd}
Pentahedral prisms which contain a \qrtet\ whose score does not
exceed $-0.52 \scorept$, or a quad
cluster whose score does not exceed $-1.04 \scorept$, or a
pentahedral cap whose score does not exceed $3.48 \scorept$, must
score less than $8 \scorept$.
\end{lem}

\begin{proof}
We begin by recalling various 
bounds from other papers in the
program.  First, by Calculation~9.1
of \cite{sp1}, $\gma(S) \le 1 \scorept$ for any
\qrtet\ $S$.
Second, recall Lemma~9.17 of \cite{sp1}, which states that if a \qrtet\
has circumradius at least $1.41$, meaning that it must be 
vor scored, then its score is less than $-1.8 \scorept$.

Together, these two results imply that the score of a \qrtet\ cannot
exceed $1 \scorept$.
Lemma~9.6 of \cite{sp1} states that if all of the \qrtets\ in a
pentahedral cap are scored by compression, the 
score of such a pentahedral cap is less than $4.52 \scorept$.
Therefore, we deduce that the score of any
pentahedral cap is less than $4.52 \scorept$.

%In the Computations section
%of this paper, we prove that the score of a quarter is nonpositive.

By Lemma~3.13 of \cite{FKC}, the score of a quad cluster is
nonpositive.  Using the scoring bound for pentahedral caps
together with this result,
we conclude that a pentahedral prism must score less than $9.04 \scorept$.
 
Therefore, if the score of a
quad cluster does not exceed $-1.04 \scorept$, the score of the
pentahedral prism to which it belongs must fall below $8 \scorept$.
Similarly, if the score of a \qrtet\ does not exceed $-0.52 \scorept$,
the score of the pentahedral cap containing it cannot exceed
$(4-0.52) \scorept$, bringing the score of the associated pentahedral
prism below $8 \scorept$, as the contribution from the other
pentahedral cap is less than $4.52 \scorept$.
In particular, if any \qrtet\ is vor scored, the score of the
pentahedral prism to which it belongs must fall below $8 \scorept$.

Likewise, if the score of a pentahedral cap does not exceed
$3.48 \scorept$, the score of the pentahedral prism with which it
is associated must fall below $8 \scorept$.
\end{proof}

We restrict our attention to pentahedral prisms not treated
by Lemma~\ref{lem:pentscorebd}.

\begin{remark}
As usual, we have chosen the constant $4.52 \scorept$ as a relatively
simple representative of many possible scoring bounds.  We could prove
a slightly tighter bound, but the resulting number would be more
cumbersome.
\end{remark}

We now present the main result of this paper.

\begin{thm}
The score of the pentahedral prism is less than $8 \scorept$.
\end{thm}

\begin{proof}
We prove linear relations on all of the 
standard clusters in the pentahedral
prism.  We combine these relations to prove the required scoring
bound.  We invoke Lemma~\ref{lem:pentscorebd} repeatedly, in
order to simplify the numerical verification of these relations.

In the section titled Computations, we establish
\[ \score + m \sol - b \le 0, \]
for all quad clusters 
not treated by Lemma~\ref{lem:pentscorebd}.
Here $\score$ denotes the
score of a quad cluster, $\sol$ denotes the spherical angle associated
with the quad cluster, $m = 0.3621$ and $b = 0.49246$.
%line1[x_]:= 0.4922197796533495 - 0.3621*x

In addition, we establish 
\[\score + m \sol + \epsilon (\dih - \frac{2 \pi}{5}) - b_c \le 0\]
for all
compression-scored \qrtets\ forming part of a pentahedral cap
not treated by Lemma~\ref{lem:pentscorebd}.
Here $\dih$ denotes the dihedral angle associated with the first edge
of the tetrahedron (that is, the edge common to the five tetrahedra in a
pentahedral cap), $\epsilon = 0.0739626$ and $b_c = 0.253095$.
%tetshift[x_]:= 0.253095 - 0.3621*x (eps = 0.0739626)

Summing the relations for the five \qrtets\ from a pentahedral cap,
we find
\[
\sum_{i=1}^5 \score_i + m \sum_{i=1}^5 \sol_i + 
\epsilon \sum_{i=1}^5 (\dih_i - \frac{2 \pi}{5}) - 5 b_c \le 0.
\]
Summing over both pentahedral caps and using the relation that
the sum of the five dihedral angles in a pentahedral cap is $2 \pi$,
\[
\sum_{i=1}^5 \dih_i = 2 \pi,
\]
we find
\[
\sum_{i=1}^{10} \score_i + m \sum_{i=1}^{10} \sol_i - 10 b_c \le 0.
\]
We represent the tetrahedra from the second
pentahedral cap by the indices $i=6 \ldots 10$.

Adding the relations for
the five quad clusters (indexed from $11$ to $15$), 
and using the fact that the sum of the solid
angles is $4 \pi$,
\[
 \sum_{i=1}^{10} \sol_i + \sum_{j=11}^{15} \sol_j = 4 \pi
\]
we find
\[
\sum_{i=1}^{10} \score_i  + \sum_{j=11}^{15} \score_j + 
4 \pi m - 5 b - 10 b_c \le 0.
\]

Therefore,
\[
\score \le 5 b + 10 b_c - 4 \pi m.
\]
The left-hand side denotes the score of the pentahedral prism.  If the right-hand
side is bounded below $8 \scorept$, we have achieved the required result.
Substituting the values of $b, b_c, m, \text{ and } \scorept$, we find that the
score of the pentahedral prism is less than $7.9997 \scorept$. %\qed
\end{proof}

\begin{remark}
The majority of the bounds which we prove in this paper are
``relaxed'', in the sense that they are $\epsilon$-away from the
ideal bound.  By decreasing $\epsilon$, we can prove tighter bounds,
at the cost of increasing the complexity of the computer verifications.

We could prove a slightly tighter bound on quad clusters (by decreasing $b$).
This would improve the bound on pentahedral prisms slightly, to 
$7.98 \scorept$, say.
\end{remark}

\section{Interval Arithmetic}
Due to the complex nature of the local density functions, it is typically
unrealistic to attempt to prove the required relations directly.  Instead,
we prove the majority of the required relations via
computer-based interval arithmetic methods.

We review the basic notions of interval arithmetic.

Suppose that the value of a function $f(x)$ lies in the interval
$[a,b]$.  Further, suppose that $g(x)$ lies in the interval $[c,d]$.
Then $f(x) + g(x)$ must lie in $[a+c,b+d]$.  While it may be the case
that we could produce better bounds than this for the function $f + g$,
these interval bounds give crude control over the behavior of the
function.  Interval arithmetic provides a mechanism for formalizing
arithmetic on these bounds.

We represent an interval $t$ as $\big[\underline{t},\overline{t}\big]$.
Then for intervals $a$ and $b$,
\[
a + b = \big[\underline{a} + \underline{b}, \overline{a} + \overline{b}\big].
\]
Likewise,
\[
a - b = \big[\underline{a} - \overline{b}, \overline{a} - \underline{b}\big].
\]
Multiplication is somewhat more complicated.  Define 
\[
C = \{\underline{a} \/ \underline{b}, \underline{a} \overline{b},
	\overline{a} \underline{b}, \overline{a} \overline{b}\}.
\]
Then
\[
a*b = [\min(C), \max(C)].
\]
We leave division as an exercise for the reader.

Similarly, we can define the operation of a monotonic function on
an interval.  For example,
\[
\arctan(a) = [\arctan(\underline{a}), \arctan(\overline{a})].
\]

Using interval arithmetic, we can produce crude bounds for polynomials
evaluated on intervals.  Likewise, we can produce crude bounds for
rational functions evaluated on intervals.  Finally, we add the
composition of monotonic functions.  This allows us to produce interval
bounds for functions such as $\sol$ and $\vor$ over \qrtets,
quarters, or quad clusters.

\section{The Method of Subdivision}
The relations on tetrahedra and quad clusters required for
the scoring bound on decomposition stars typically have
the form
$g(y)\le 0$ for $y \in I$, where $I$ is a product of closed
intervals.  As $g$ is usually continuous, the existence of a
maximum is trivial.  However, bounds on the behavior of $g$
over all of $I$ computed directly via interval arithmetic
are generally poor.

We define a {\em cell} to be a product of closed intervals.
By subdividing $I$ into
sufficiently small cells, the quality of the computed bounds on
each cell usually
improves enough to prove the relation for each cell, and hence
for the original domain $I$.

If in fact $g(y) \le c < 0$, this approach works very well.  However,
if the bound is tight at a point $y_0$, i.e., $g(y_0)=0$, then pure
subdivision will usually fail, since the computed upper bound on $g$ over
any cell containing $y_0$ will typically be positive.

If $y_0$ is not an interior maximum, 
we turn to the partial derivatives of $g$.  If
we can show that the partials of $g$ on a small cell containing
$y_0$ have fixed sign (bounded away from zero), then the
maximum value of $g$ on that cell is easily computed.  
It is typically the case that a cell
must be very small before we can determine the sign of the
partials via interval arithmetic bounds.

\section{Dimension Reduction}
The relations on tetrahedra required for
the scoring bound on decomposition stars are typically
six-dimensional.  For a quad cluster, they can be even
higher-dimensional.  For high-dimensional relations,
the method of subdivision becomes very expensive,
computationally speaking.

We define a simplification which reduces the dimension
of the required computations.  This simplification
therefore reduces the computational expense of the
verification of a relation.

We refer to this simplification as
dimension-reduction.  The details of the proof of the
simplification vary depending on 
whether the scoring is by compression, vor analytic, or
Voronoi.

Introduced in \cite{sp1} for compression scoring, this
argument states that moving a vertex $v_i$ along the
edge $(0,v_i)$ toward the origin $0$
holds the solid angle fixed, while
increasing the score of the tetrahedron.  See Figure~\ref{fig:tet}.
Since the reduction may be performed until an edge-length constraint
is met, this argument reduces the number of free
parameters for the verification, thus reducing the dimension and
complexity of the verification of a relation.

The validity of the same reduction for vor analytic-scored
tetrahedra is obvious if the
tip of the Voronoi cell does not protrude.
If the tip does protrude, we must use the analytic continuation
for the Voronoi volume.  In this case, the validity of
the reduction is not obvious.

We provide a sketch of an analytic proof that this reduction increases
the analytic Voronoi score of a tetrahedron.
The geometric constraint of moving a vertex along an edge
can easily be stated analytically in terms of the original edge
lengths, $(y_1, y_2, y_3, y_4, y_5, y_6)$.
This action depends on a single parameter, the
distance of the vertex $v_1$ from the origin, which we call $t$.
The new edge lengths are given by
\[
( t, y_2, y_3,
  y_4,\sqrt{{t^2} + {{{y_3}}^2} - 
   {\frac{t\,\left( {{{y_1}}^2} + 
         {{{y_3}}^2} - {{{y_5}}^2}
         \right) }{{y_1}}}},
  \sqrt{{t^2} + {{{y_2}}^2} - 
   {\frac{t\,\left( {{{y_1}}^2} + 
         {{{y_2}}^2} - {{{y_6}}^2}
         \right) }{{y_1}}}} ).
\]

Recall from \cite{sp1} that the formula for the analytic Voronoi
volume is a rational function of $\chi$, $u$, $\sqrt{\Delta}$, and
$x_i$, where $x_i = y_i^2$.  Further recall that $\chi$, $u$, and
$\Delta$ are all polynomial functions in $x_i$.

Substituting the computed edge lengths in the formula for the analytic
Voronoi volume, taking the partial derivative with
respect to $t$, replacing $t$ with $y_1$, multiplying by the positive
term
\[
8 \sqrt{\Delta} u(x_1, x_3, x_5) u(x_1, x_2, x_6)/y_1,
\]
and then simplifying, we end up with a large homogeneous
polynomial in $x_i$ of degree 6, which is too ugly to exhibit here (having
91 terms).

Evaluating this polynomial over all possible
\qrtets\ and quarters, we find that it is positive.

Therefore the
volume is increasing in $t$, so to increase the score, we should
push the vertex in along the edge.
The verification of the sign of the polynomial is found in
Calculation~\ref{d:dimred}.

The validity of a similar reduction argument for Voronoi scoring
of a quad cluster
is obvious, since the Voronoi volume is increasing in $t$.

%\begin{figure}
%\includegraphics{tet.eps}
%%\centerline{ \psfig{file=tet.eps} }
%\caption{Tetrahedron with distinguished vertex and labeled edges.}
%\label{fig:tet}
%\end{figure}

\begin{remark}
If the computational effort to prove the relations we require
for the scoring bound on decomposition stars were not so
extreme, we could dispense with the complications associated with
dimension-reduction.
\end{remark}

\section{Computations}
We derive the auxiliary bounds necessary to prove the scoring bound.  We
separate the bounds into sections, depending on the case to
which each bound is to be applied.

\begin{remark}
We go to great lengths to reduce the complexity of the relations
which we are required to verify.  If we had sufficient computer
resources, we could dispense with many of the
methods which we use to reduce the complexity of the calculations.
This would simplify the computations significantly.
However, the majority of the relations which we wish to prove
are complex enough that a direct approach rapidly overwhelms our
available resources.
\end{remark}

\subsection{Quarters}
We are required to prove that the score of a quarter is
nonpositive.

Quarters are either flat or upright.  Flat quarters are scored using
either compression or vor analytic scoring.  Upright quarters 
occurring in a quad cluster are
scored using either compression or averaged vor analytic scoring.
Recall the scoring scheme from \ref{scoring:quarter}
in the Introduction. 

Since compression does not depend on a distinguished vertex, we need
only consider the flat case.  Calculation~\ref{q:gma} shows that the
compression score of any quarter is nonpositive. 

Averaged vor analytic scoring, applied to an upright quarter, is
similarly simple.  Calculation~\ref{q:octavor} shows that the averaged 
vor analytic score of any upright quarter is nonpositive.

\begin{remark}
In these first two calculations, we have
simplified the calculations by proving a 
stronger result than that which is strictly necessary.  
In the last case, we
are required to fully invoke the scoring scheme.  This complicates the
analysis.
\end{remark}

The only case remaining is that of a flat quarter, scored by 
vor analytic.  This case is trickier, since it is not true that
vor analytic is nonpositive on any flat quarter.  Therefore, we
must assume that one of the two faces adjacent to the diagonal
has circumradius not less than $\sqrt{2}$.  This constraint
complicates the verifications.

We first consider a small cell containing the
edge lengths $(2,2,2,2,2,2\sqrt{2})$.  We call a cell containing
these edge lengths a
{\em corner} cell.  We prove in Calculation~\ref{q:partials} that for a
sufficiently small corner cell, the $y_1$ through $y_5$ partials
(the edge lengths corresponding to the short edges) are negative.

To find the maximum score, we therefore decrease the short edges
as much as possible, while not violating the face constraint.
On such a cell, the face constraint of one of the edges is tight,
so we may assume that $\eta(y_1,y_2,y_6)=\sqrt{2}$ and $y_3=2$,
$y_4=2$, and $y_5=2$ or else $\eta(y_4,y_5,y_6)=\sqrt{2}$ and
$y_1=2$, $y_2=2$, and $y_3=2$,
where $\eta(y_1,y_2,y_3)$ is the circumradius of a face with edge
lengths $(y_1,y_2,y_3)$.

We represent the desired relation
as $sc + \alpha (\eta^2 - 2) \le 0$ for $\alpha = 0.125$.
%define ALPHA			0.125
Calculation~\ref{q:corner} verifies this relation.

Second, we prove the desired relation off of the corner cell.  We
subdivide this verification into Calculation~\ref{q:off:dimred} 
taking advantage of
dimension reduction and partial derivative information, 
and Calculation~\ref{q:off:bdry} which only considers the boundary 
(where the face constraint is tight).

\subsection{Quasi-regular Tetrahedra}
%Otherwise known as the pentahedral cap.

There are three verifications required to prove the desired
relation on \qrtets.  First, we prove a relation between dihedral
angle and score.  We then show that if the dihedral angle of a
tetrahedron in a pentahedral cap exceeds a certain bound, the score
of the pentahedral cap must fall below $3.48 \scorept$, 
%hence such arrangements may be discarded.
falling into the purview of Lemma~\ref{lem:pentscorebd}.
We call such a bound a
{\em dihedral cutoff}.  This cutoff then allows us to prove the final bound.

In the following discussion, $\dih$ refers to the dihedral angle
associated with the first edge of a \qrtet, $\score$ refers to
the compression score of the tetrahedron, 
and $\sol$ refers to the solid
angle at the distinguished vertex.  We restrict
our attention to \qrtets\ whose score exceeds $-0.52 \scorept$, as
per Lemma~\ref{lem:pentscorebd}.

The first relation has the form $\score \le a_1 \dih - a_2$,
where $a_1 = 0.3860658808124052$ and $a_2 = 0.4198577862$.
Calculation~\ref{qr:dihrel}
provides the verification of this relation.

%smside2[tt_]:= -0.4198577860315257 + 0.3860658808124052*tt

Applying the
relation to four \qrtets\ forming part of a pentahedral
prism, we find
\[
\sum_{i=1}^4 \score_i \le a_1 \sum_{i=1}^4 \dih_i - 4 a_2.
\]
Applying the relation
\[
\dih_5 = 2 \pi - \sum_{i=1}^4 \dih_i
\]
and adding $\score_5$ to both sides of the first relation,
we find
\[
\sum_{i=1}^5 \score_i \le \score_5 + a_1 (2 \pi - \dih_5) - 4 a_2 .
\]
The left-hand side represents the score of the pentahedral cap.
If the right-hand side does not exceed $3.48 \scorept$, we
can remove the arrangement from consideration, since it pulls the score of
the associated pentahedral prism below $8 \scorept$, as per
Lemma~\ref{lem:pentscorebd}.

%d_0 = 1.4674

We assert that if $\dih \ge d_0$, where $d_0 = 1.4674$,
the right-hand side
\[
\score_5  + a_1 (2 \pi - \dih_5) - 4 a_2
\]
does not exceed $3.48 \scorept$.  Rewriting this statement,
we prove that
$\dih \ge d_0$ implies
\[
\score - a_1 \dih \le 3.48 \scorept - 2 \pi a_1 + 4 a_2,
\]
which is verified in Calculation~\ref{qr:dihcut}.
Hence we may restrict our attention to \qrtets\ whose
dihedral angle does not exceed the dihedral cutoff $d_0$.

Using the dihedral cutoff, we establish the final relation,
\[\score + m \sol + \epsilon (\dih - \frac{2 \pi}{5}) - b_c \le 0.\]
Calculation~\ref{qr:pentacap} provides the verification.

\subsection{Flat Quad Clusters}
Flat quad clusters are composed of two flat quarters, whose common face
includes the long edge.
We prove the relation $\score \le -m \sol + b$
on quad clusters whose score exceeds $-1.04 \scorept$, again
invoking Lemma~\ref{lem:pentscorebd}.
We arrive at this relation
for flat quad clusters
by proving the relation $\score \le -m \sol + b/2$ on flat quarters.
Here $\score$ refers to vor or compression scoring, whichever is appropriate
for the quarter.

We restrict our attention to flat quarters whose score exceeds
$-1.04 \scorept$, recalling the fact that the score of
flat quarters is non-positive.
Adding the relation for each flat quarter, we arrive at the desired
bound for flat quad clusters.

In the following discussion, we label the diagonal of a flat
quarter $y_6$.

Flat quarters may be scored using either compression or vor
scoring.  We treat each case separately.

First, suppose that we wish to prove the bound for compression scored
quarters.  This means that the circumradii of the two faces adjacent
to the long diagonal do not exceed $\sqrt{2}$.  We subdivide the
verification into 
Calculation~\ref{flat:gma:dimred}, 
a computation where we apply dimension-reduction
and partial derivative information, and 
Calculation~\ref{flat:gma:bdry}, a boundary verification, where
we restrict our attention to cells which lie on the boundary between
compression and vor scoring.

%Calculation~\ref{flat:gma:dimred} provides the verification 
%of the dimension-reduction
%case.  Calculation~\ref{flat:gma:bdry} provides the 
%verification of the boundary case.

Second, we treat the vor-scoring case.  In this case we prove the
bound for vor-scored quarters.  This means that at least one of
the circumradii of the two faces adjacent to the long diagonal
is at least $\sqrt{2}$.  This verification is somewhat more complex
than the compression case.  We subdivide the verification into
\begin{enumerate}
\item 
Verification that the first three partials are negative on a
small cell containing the corner (Calculation~\ref{flat:partials}).
\item 
Verification of the bound on that small cell containing the corner,
using the property that the first three partials are negative
(Calculation~\ref{flat:corner}).
\item 
A computation where we apply dimension-reduction and partial
derivative reduction, omitting the corner cell 
(Calculation~\ref{flat:vor:dimred}).
\item 
A boundary verification, where we restrict our attention to
cells which lie on the boundary between compression and vor
scoring, again omitting the corner cell
(Calculation~\ref{flat:vor:bdry}).
\end{enumerate}

%These verifications are contained in Calculations (ref.).

\subsection{Octahedra}
%Four upright quarters, arrayed around their common long edge so
%that each face containing the common edge is shared by two quarters,
%form an {\em octahedron}, another type of quad cluster.

Recall that octahedra, a type of quad cluster,
are composed of four upright 
quarters arrayed around their common long edge 
(known as the {\em diagonal}) so that each
face containing the common edge is shared by two quarters.

We are required to prove a relation of the form 
\[ \score + m \sol - b \le 0, \]
where $\score$ denotes the score of an octahedron, $\sol$ denotes
the solid angle associated with the distinguished vertex, and $m$
and $b$ are positive constants.  By Lemma~\ref{lem:pentscorebd}, 
we restrict our attention to
octahedra whose score exceeds $-1.04 \scorept$.

Our treatment of octahedra, as usual, is comprised of a number of auxiliary
computations.  We prove bounds on upright quarters which are part of an
octahedron, and then combine these bounds to deduce the required bound on
octahedra in general.

\begin{remark}
Due to the complex nature of octahedra, we are
required to consider a number of sub-cases, which we hope will not
fatigue the reader.  These cases are partitioned according to the length
of the diagonal and the scoring system applied to the upright quarters.

Using a dihedral summation argument, we will eliminate octahedra
whose diagonal lies in the range $[2.51, 2.716]$.

Next, we will treat
the case where the diagonal lies in the range $[2.716, 2\sqrt{2}]$.
Using a dihedral correction term, 
we will prove the bound for octahedra which are completely
compression-scored, and octahedra which are completely vor-scored.

The remaining cases will consist of
octahedra which contain either two or three
vor-scored quarters.  (Since a quarter is vor-scored if one
of the faces containing the diagonal has circumradius $\sqrt{2}$ or
greater, it is not possible for an octahedron to contain only one
vor-scored quarter.)  We treat these cases using an additional
correction term.

The details follow.
\end{remark}

In all computations involving octahedra, 
we label the diagonal $y_1$.

In order to simplify the computations, 
we first prove an auxiliary cutoff bound.
This first bound reduces
the size of the cell over which we must conduct our search,
as per Lemma~\ref{lem:pentscorebd}.

Calculation~\ref{octa:peel} shows that if an upright quarter
contains an edge numbered
2, 3, 5, or 6 whose length is not less than $2.2$, its score is less than
$-0.52 \scorept$.

Since such an edge is shared by another upright quarter
in the same octahedron, the score of the 
associated octahedra must fall below $-1.04 \scorept$.

We restrict our search accordingly.

%Octahedral cutoff:  show that if the diagonal of an
%octahedron lies in the range [2.51,2.71], the score of
%the octahedron is less than -1.04 pt.  Hence, we can
%discard such arrangements.
%Specifically, we verify the following bound to within
%a tolerance of (fill_in_the_blank):
%sc < m*dih + b,
%where
%Newer version:  cutoff = 2.716,
%{m, b} = {-0.1533667634670977, 0.2264803995076098}
%Required tolerance (for the new bound):  3.0e-5

In the first case, we assume that the diagonal lies in the range
$[2.51, 2.716]$.  In Calculation~\ref{octa:cut}, 
we prove a bound of the form
\[
\score + c \dih \le d
\]
on upright quarters, where $c = 0.1533667634670977$, and
$d = 0.2265$.  Adding the bound for four quarters forming
an octahedron, we find
\[
\sum_{i=1}^4 \score_i + c \sum_{i=1}^4 \dih_i \le 4 d.
\]
Using the fact that the sum of the dihedral angles is $2 \pi$,
we find that
\[
\score \le -2 \pi c + 4 d.
\]
A computation involving the constants $c$ and $d$ shows that the
score is less than $-1.04 \scorept$.  Again invoking
Lemma~\ref{lem:pentscorebd}, we need only
consider octahedra whose diagonal lies in the range 
$[2.716, 2\sqrt{2}]$.

%Next, we consider octahedra which contain (at least two)
%vor-scored quarters whose diagonal falls in the range
%$[2.716, 2.81]$.  We prove a cutoff bound, showing that
%a vor-scored quarter with such a diagonal must score
%less than $-0.52 \scorept$.  Since there must be at least
%two such quarters if there are any, the score of the
%octahedron with which they are associated must fall below
%$-1.04 \scorept$.  We therefore discard such arrangements,
%restricting our attention to diagonals in the range

Using this assumption, we prove bounds of the form
\begin{equation}
%\[
\score + m \sol + \alpha \dih  \le \frac{b}{4} + 
	\alpha \frac{\pi}{2}
%\]
\label{octa:dih}
\end{equation}
and
\begin{equation}
%\[
\score + m \sol + \alpha \dih + \beta x_1 \le \frac{b}{4} + 
	\alpha \frac{\pi}{2} + 8 \beta,
%\]
\label{octa:edge}
\end{equation}
where $\dih$ refers to the dihedral angle associated with the diagonal,
$\score$ refers to the scoring scheme appropriate for a
particular upright quarter, and $x_1$ refers to the square of the length of
the diagonal.  We
choose $\alpha$ and $\beta$ according to the scoring scheme.

\begin{remark}
Appropriate values for the correction terms involving 
$\alpha$ and $\beta$ were determined by experimentation.
\end{remark}

%For the first inequality, we choose $\alpha = 0.14$ for compression-scored
%quarters, and $\alpha = 0.054$ for vor-scored quarters.

%We prove Relation~\ref{octa:dih} for compression-scored quarters, choosing
%$\alpha = 0.14$.

Choosing $\alpha = 0.14$, we prove (\ref{octa:dih}) for compression-scored
quarters with diagonal in the interval $[2.716, 2\sqrt{2}]$
(Calculation~\ref{octa:gma:dih}).  Using
the same $\alpha$, we prove (\ref{octa:dih}) for vor-scored
quarters with diagonal in the range $[2.716, 2.81]$ (Calculation~\ref{octa:vor:dih}).

Choosing $\alpha = 0.054$, $\beta = 0.00455$, we prove (\ref{octa:edge}) for
compression-scored quarters with diagonal in $[2.81, 2\sqrt{2}]$
(Calculation~\ref{octa:gma:corr}).  Choosing
the same $\alpha$, but $\beta = -0.00455$, we prove (\ref{octa:edge}) for
vor-scored quarters with diagonal in $[2.81, 2\sqrt{2}]$
(Calculation~\ref{octa:vor:corr}).

%For the second inequality, we choose $\alpha = 0.054$, and we choose
%$\beta = 0.00455$ for compression-scored quarters, and $\beta = -0.00455$
%for vor-scored quarters.

Note that for vor-scored quarters,
the first inequality is a relaxation of the second, since $\beta$ is
negative.

The verification of each of these inequalities involves
a computation where we apply dimension-reduction
and partial derivative information, and 
a boundary verification, where
we restrict our attention to cells which lie on the boundary between
compression and vor analytic scoring.
Note that the
dimension-reduction step for relation~(\ref{octa:edge}) is complicated 
by the presence of the $\beta x_1$ term.

Summing the first inequality over an octahedron, we find
\[
\sum_{i=1}^4 \score_i + m \sum_{i=1}^4 \sol_i + \alpha \sum_{i=1}^4 \dih_i
\le b + 2 \alpha \pi.
\]
Using the fact that the dihedral angles sum to $2 \pi$, we find
\[
\score  + m \sol \le b,
\]
so octahedra with diagonals in the range $[2.716,2.81]$ 
satisfy the requisite bound.

Summing the first inequality over a consistently scored octahedron
(either all compression or all vor) with diagonal in the
range $[2.81, 2\sqrt{2}]$, we again arrive at the desired bound.

The remaining cases involve octahedra which contain both compression
and vor-scored quarters, and whose diagonals lie in the range
$[2.81, 2\sqrt{2}]$.  For this case, we use the second inequality.

The summation involving the second inequality is identical to the
first, save for
the presence of the $\beta$ terms.  If there are two vor-scored
quarters and two compression-scored quarters, the beta terms cancel, 
giving the relation as before.

If there are three vor-scored quarters and one compression-scored
quarter, we note that the same relation for vor-scored quarters
holds if we replace $\beta$ by $\beta/3$ (since we have now relaxed
the bound).  Summing the inequalities, the term involving $\beta$
vanishes again, leaving the desired inequality.

\subsection{Pure Voronoi quad clusters}
The next class of quad clusters which we treat are the pure
Voronoi quad clusters.  We will define a truncation operation
on these quad clusters.  Truncation will simplify the geometry of
the quad clusters, and will provide a convenient scoring bound.
We will then divide our treatment of pure Voronoi quad clusters
into two cases in order to simplify the analysis and numerical
verifications as much as possible.

\begin{figure}
\includegraphics{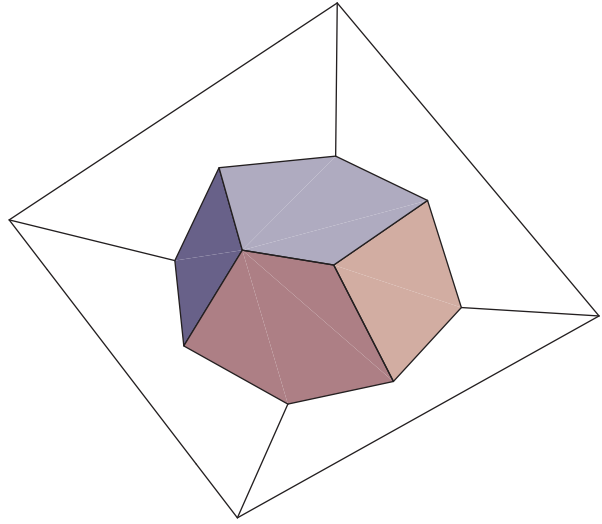}
%\centerline{ \psfig{file=vor.eps} }
\caption{A pure Voronoi quad cluster.}
\label{fig:vor}
\end{figure}

Recall from the classification of quad clusters (\ref{quad:classification})
that a pure Voronoi quad cluster consists of the intersection of a
$V$-cell at the origin with the cone at the origin over a quadrilateral
standard region.  We refer to the restriction of the $V$-cell to the
cone over the quadrilateral as either the $V$-cell 
or the Voronoi cell of the quad cluster.
Figure~\ref{fig:vor} describes the geometry of a simple
$V$-cell.

In addition, recall that a vertex lying in the cone over a pure
Voronoi quad cluster must have height greater than $2\sqrt{2}$.
Such vertices can significantly complicate the geometry of the $V$-cell,
affecting its shape and volume.

We remove the effect of vertices lying above a pure Voronoi quad cluster
by removing all points from the $V$-cell which have height greater than
$\sqrt{2}$.  We call this operation {\em truncation at} $\sqrt{2}$.
Truncation decreases the volume of the quad cluster.  This decrease in
volume increases the score of the quad cluster, bringing it closer to
the proposed bound.

We refer to truncated pure Voronoi quad clusters as {\em truncated} quad
clusters.

We define a scoring operation on pure Voronoi quad clusters which we
call {\em truncated Voronoi} scoring.  This operation consists of
truncation at $\sqrt{2}$, followed by the usual Voronoi scoring.

Each diagonal across the face of a cluster must have length
greater than $2\sqrt{2}$, otherwise we could form two flat quarters,
contradicting the decomposition.  We choose the shorter of the two
possible diagonals, and will consider that diagonal in the analysis
which follows.

We decompose the cluster into two tetrahedrons along the chosen
diagonal.  The face dividing the tetrahedrons is either acute or
it is obtuse.  We treat each case separately.

We must prove 
\[ 
\score + m \sol - b \le 0, 
\]
where $\score$ denotes the score of a pure Voronoi quad cluster, $\sol$ denotes
the solid angle associated with the distinguished vertex, and $m$
and $b$ are positive constants.
We call this relation a {\em bound} on the solid angle and score of a quad
cluster.
Invoking Lemma~\ref{lem:pentscorebd}, 
we restrict our attention to
quad clusters whose score exceeds $-1.04 \scorept$.

\subsubsection{The acute case}
If the separating face is acute, we prove
\[ \score + m \sol - b/2 \le 0 \]
for each half
independently, and deduce the desired bound by adding the bounds
for each half.  Since the score of each half is non-positive (by
the arguments of Lemma~3.13 of \cite{FKC}),
we may restrict our attention to halves whose score exceeds
$-1.04 \scorept$, by Lemma~\ref{lem:pentscorebd}.

The required bound has
the property that
\[ m \sol_0 - b/2 \le 0, \]
where $\sol_0$ denotes the solid angle of the tetrahedron
$(2,2,2,2,2,2\sqrt{2})$.  If $\sol < \sol_0$, then
$m \sol < m \sol_0$, hence
\[m \sol -b/2 < m \sol_0 - b/2 \le 0,\] 
and
\[sc + m \sol -b/2 < sc \le 0,\] 
so the bound follows.  We therefore
may restrict our attention to halves whose solid angle is at least
$\sol_0$.  In addition, we restrict our attention to halves for
which the dividing face is acute.

The required verifications for each half of an acute quad cluster
are somewhat difficult to achieve
directly, so we subdivide into a number of different cases 
in an attempt to reduce the complexity of
the calculations.  First, we show that the bound holds for
all halves whose diagonal is at least $2.84$ 
(Calculation~\ref{acute:cut}).  
Using this information, we then prove the
bound everywhere but in a small corner cell
(Calculation~\ref{acute:vor}).  We then restrict our
attention to the small corner cell (Calculation~\ref{acute:corner}).  
These computations involve
the use of partial derivative information, and include the required
boundary computations.

\subsubsection{The obtuse case}
If the separating face is obtuse, the analysis becomes significantly
harder.  It is no longer possible to prove the desired bound on each
half independently.  The dimension of the full bound, even using
the usual dimension-reduction techniques, is too high to make the
verification tractable numerically.  Therefore we adopt a different
approach.

Using the dimension-reduction technique, we push each vertex along its
edge until the distance from each vertex to the origin is $2$.  
We call the resulting quad cluster a {\em squashed} cluster.  Observe
that the solid angle of the cluster is unchanged, while the volume of
the Voronoi cell has decreased, thereby increasing the score of the
cluster.

Since
the central face is still obtuse, the length of the diagonal after
this perturbation must still exceed $2\sqrt{2}$.  Note, however,
that the other edge lengths in the quad cluster can be as small as
$4/2.51$.

\begin{figure}
\includegraphics{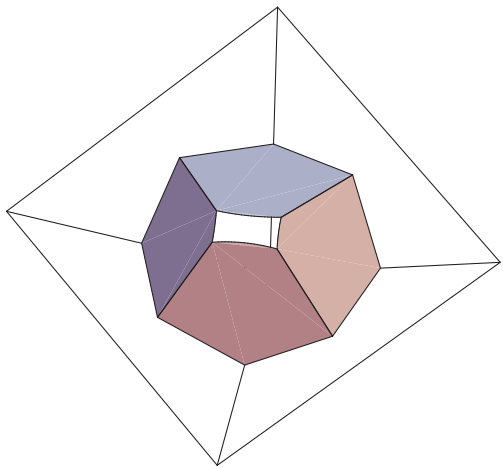}
%\centerline{ \psfig{file=vor3.eps} }
\caption{A typical truncated quad cluster.}
\label{fig:vor3}
\end{figure}

\begin{figure}
\includegraphics{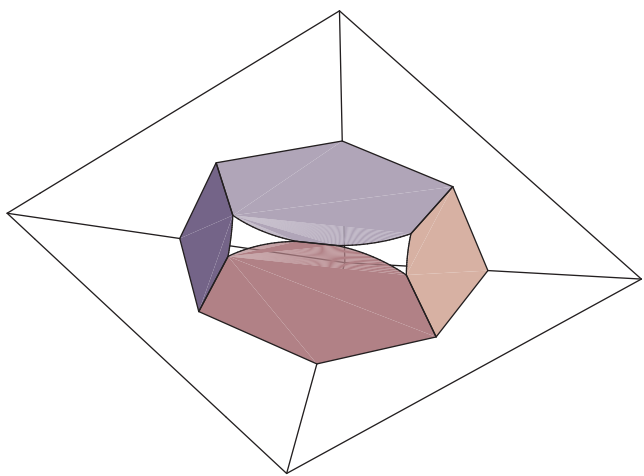}
%\centerline{ \psfig{file=vor2.eps} }
\caption{An impossible arrangement.}
\label{fig:vor2}
\end{figure}

The geometry of the $V$-cell of a squashed cluster, assuming that
there is no truncation from vertices of the packing lying above the quad
cluster, is that of Figure~\ref{fig:vor}.  When the $V$-cell is
truncated at $\sqrt{2}$ from the origin, two potential arrangements
arise.  In the first arrangement, the truncated region is connected,
as in Figure~\ref{fig:vor3}.  In second potential arrangement, the
truncated region is formed of two disjoint pieces, 
as in Figure~\ref{fig:vor2}.  

We will conclude that the second, disjoint
case cannot arise for squashed quad clusters.  
Suppose that it could.  Pick an untruncated point
along the central ridge of the $V$-cell (see Figure~\ref{fig:vor2}).  
The distance of this point
from the origin is then less than $\sqrt{2}$, but due to its location
on the central ridge, it is equidistant from the two nearest vertices
and the origin.  This implies that the circumradius of the resulting
triangle must be less than $\sqrt{2}$, which contradicts the fact
that the diagonals have length at least $2\sqrt{2}$.

\subsubsection{A geometric argument}
We introduce a simplification which will reduce the complexity of
the obtuse case.  This simplification will consist of a perturbation of
the upper edge lengths of a squashed quad cluster.  This perturbation
will increase the score while holding the
solid angle of the quad cluster fixed.

This simplification is based on a geometric decomposition of the
truncated Voronoi cell.  We will describe the decomposition, and
then describe a construction which will ultimately simplify the
analysis.

\begin{remark}
This simplification replaces the less-rigorous construction known as the
``pizza'' argument found in previous versions of this paper.
\end{remark}

While our arguments will extend to treat a general
squashed and truncated Voronoi cell associated with a general
standard cluster, we restrict our attention to truncated
Voronoi cells associated with quad clusters.

To begin, we consider the decomposition of
a truncated Voronoi cell into its fundamental components.
A truncated Voronoi cell is formed of three elements:  a
central spherical section (formed by the truncation), {\em wedges} of
a right circular cone, and tetrahedrons called {\em Rogers simplices}.  

\begin{figure}
\includegraphics{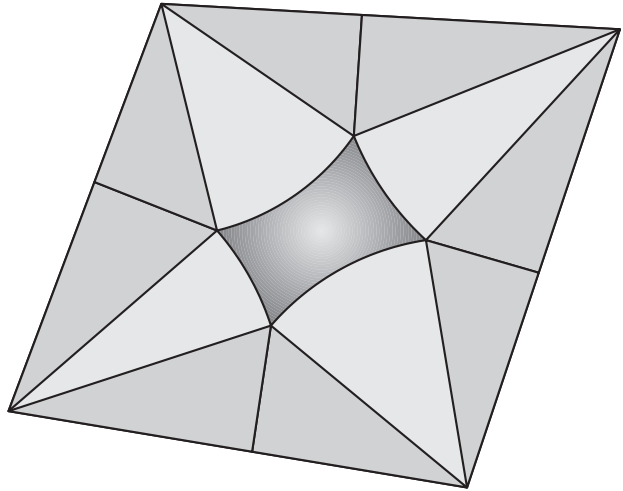}
%\centerline{ \psfig{file=quad2.eps} }
\caption{A representation of a truncated quad cluster.}
\label{fig:quad2}
\end{figure}

We choose a representation of a truncated quad cluster composed of 
the radial projection of each element to a plane passing close to
the four corners of the quad cluster.
This decomposition
is represented in Figure~\ref{fig:quad2}.

\subsubsection{Rogers simplices}
We now consider the geometry of the Rogers simplices.

\begin{figure}
\includegraphics{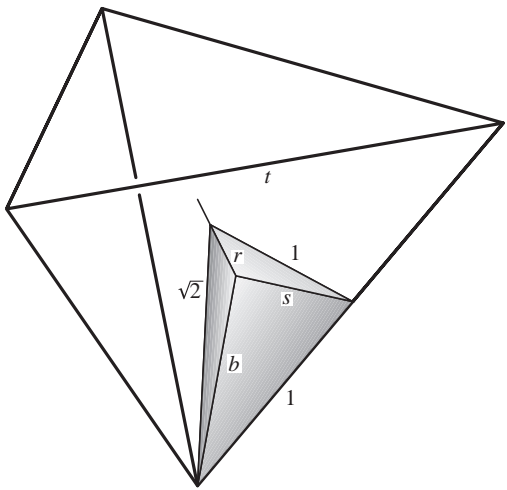}
%\centerline{ \psfig{file=quad3.eps} }
\caption{Detail of truncated Voronoi decomposition.}
\label{fig:quad3}
\end{figure}

Consider a face with edge lengths $(2,2,t)$ 
associated with a side of
a truncated quad cluster.
Let $b$ represent the circumradius of the
face, and let $r$ represent the orthogonal extension of a
Rogers simplex from the face, as in Figure~\ref{fig:quad3}.

\begin{figure}
\includegraphics{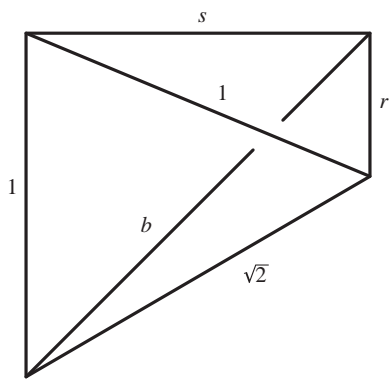}
%\centerline{ \psfig{file=quad4.eps} }
\caption{Detail of Rogers simplex.}
\label{fig:quad4}
\end{figure}

Then
\[
b = \frac{4}{\sqrt{16-t^2}}
\]
\[
r = \sqrt{2-b^2} = \sqrt{\frac{16-2t^2}{16-t^2}},
\]
and
\[
s = \sqrt{b^2-1} = \frac{t}{\sqrt{16-t^2}}.
\]
See Figure~\ref{fig:quad4}.

\subsubsection{The geometric construction}
We now present the geometric construction which will imply the
simplfication.

\begin{figure}
\includegraphics{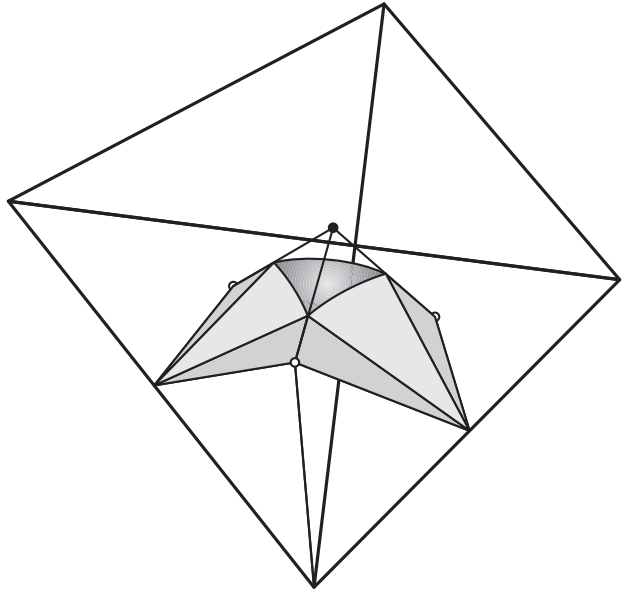}
%\centerline{ \psfig{file=decomp1.eps} }
\caption{Decomposition of a truncated Voronoi cell.}
\label{fig:decomp1}
\end{figure}

We represent the geometry of the truncated Voronoi cell associated with
one half of a quad cluster in Figure~\ref{fig:decomp1}.

\begin{figure}
\includegraphics{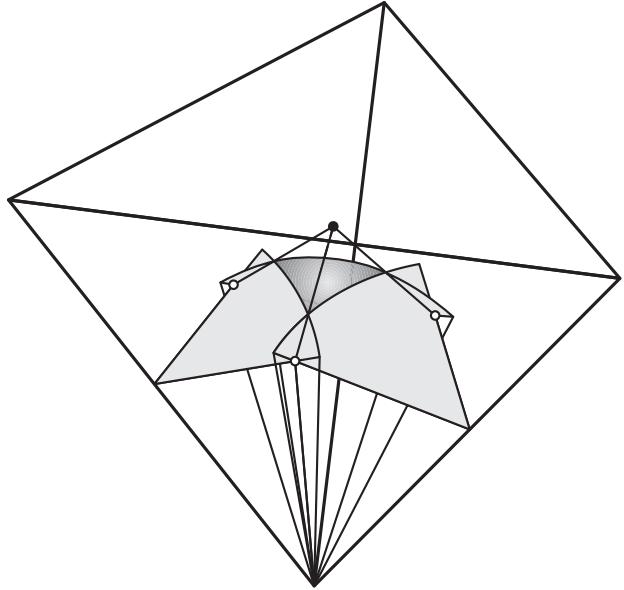}
%\centerline{ \psfig{file=decomp2.eps} }
\caption{Wedges extended to include the Rogers simplices.}
\label{fig:decomp2}
\end{figure}

We can simplify the representation by
extending the wedges to enclose the Rogers simplices.
See Figure~\ref{fig:decomp2}.
This process adds an extra volume term.

\begin{figure}
\includegraphics{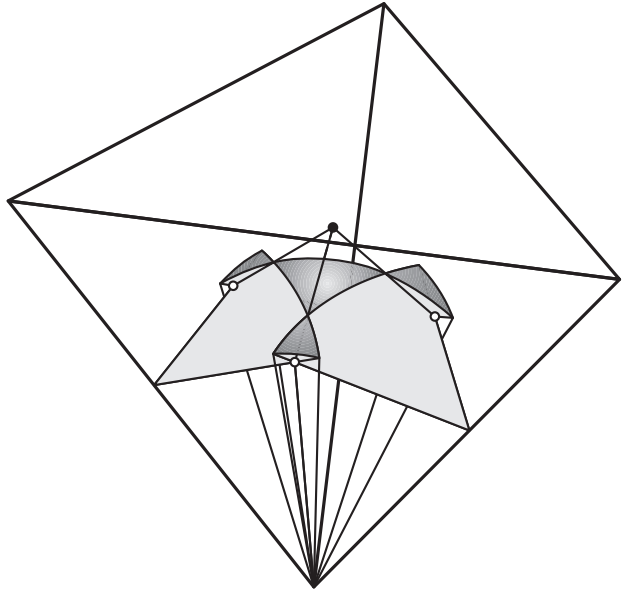}
%\centerline{ \psfig{file=decomp3.eps} }
\caption{Decomposition with flutes.}
\label{fig:decomp3}
\end{figure}

\begin{figure}
\includegraphics{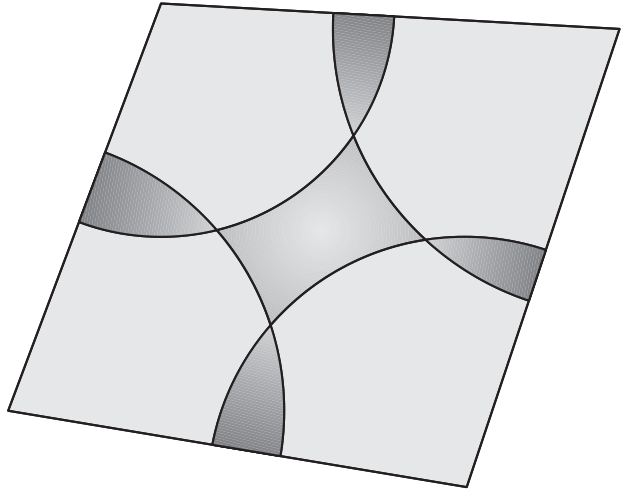}
%\centerline{ \psfig{file=planar3.eps} }
\caption{Planar representation with flutes.}
\label{fig:planar3}
\end{figure}

The overlap between the wedges is slightly complicated.  We
simplify the overlap as follows.  Take the cone over the overlap.
Intersect it with a ball of radius $\sqrt{2}$ at the origin.
We call the spherical sections produced by this construction
{\em flutes}.  This construction is represented in 
Figure~\ref{fig:decomp3}.  Figure~\ref{fig:planar3}
is a planar representation of this construction.

\begin{figure}
\includegraphics{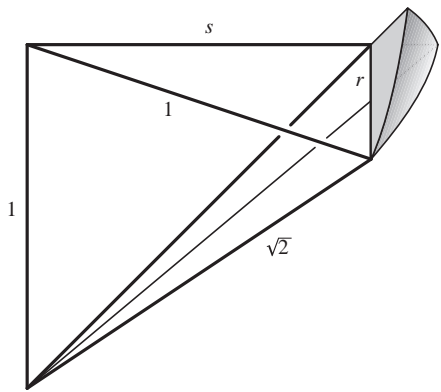}
%\centerline{ \psfig{file=quad8.eps} }
\caption{Detail of Rogers simplex with quoin.}
\label{fig:quad8}
\end{figure}

To form each flute, we have added two extra pieces of
volume (per flute) to our construction.  
We call these pieces {\em quoins}.  We
attach each quoin to a Rogers simplex.  
See Figure~\ref{fig:quad8}.

\subsubsection{A solid angle invariant}
We now require some notation for the volumes which enter
into this construction.
Let $c$ denote the volume of the central spherical angle.
Let $r$ denote the volume of the Rogers simplices.
Let $w$ denote the volume of the wedges.
Let $w'$ denote the volume of the extended wedges.
Let $q$ denote the volume of the quoins.
Let $f$ denote the volume of the flutes.
Finally, let $v$ denote the volume of the truncated Voronoi cell.
By the original decomposition,
\[
v = c + r + w.
\]
By our construction,
\[
v = c + w' + q - f.
\]

Recall that the
solid angle $s$ of the quad cluster is the sum of the 
dihedral angles minus $2\pi$.
The dihedral angles to which we refer are those 
associated with the edges between each corner
of the quad cluster and the origin.

Our perturbation will hold the solid angle $s$ of the quad cluster
fixed.  Therefore, 
the sum of the dihedral angles must also be fixed.
This fixes $w'$.

Take the cone over each extended wedge and intersect it with
a ball of radius $\sqrt{2}$ centered at the origin.  
Let $t$ denote the sum of these volumes.  Since the sum
of the dihedral angles is fixed, $t$ is also fixed.

Further, note that
\[
\frac{2\sqrt{2}}{3} s = c + t - f.
\]
This relation implies that $c-f$ is fixed.  Combining this with
the previous relations, we find that if we hold the solid angle
fixed, the volume of the truncated Voronoi cell depends only
on $q$, the volume of the quoins.

\subsubsection{The quoin}
We now develop a formula for the volume of a quoin.  We first
introduce some details of the Rogers simplices.

Consider a face $(2,2,t)$ of a truncated quad cluster.
Two Rogers simplices are associated with this face, as
suggested in Figure~\ref{fig:quad3}.  
Observe that the volume of the quoin associated with one
of these Rogers simplices is increasing
in $r = \sqrt{\frac{16-2t^2}{16-t^2}}$.  
Next, observe that $r$ is in turn decreasing in
$t$.
Therefore increasing $t$ decreases the volume 
of the quad cluster, if we hold the solid angle fixed (by varying
the length of another edge of the quad cluster).

Each half of a quad cluster has two variable edge lengths (not counting the
shared diagonal).  We label the variable
edge lengths of one half of the quad cluster $y_1$ and $y_2$.
We label the length of the diagonal $d$.
Holding the solid
angle fixed, we may perturb one half by shrinking the larger
and increasing the shorter length.  
We wish to establish that
increasing the short length reduces the volume of the 
truncated Voronoi cell
more than decreasing the longer length increases the volume.

\subsubsection{The volume of a quoin}
To achieve our reduction, we establish a formula for the volume of a
quoin.

We then verify that the volume of the quoin associated with the
shorter edge is decreasing
faster under this perturbation than the volume 
of the quoin associated with the
longer edge is increasing.

In other words, we wish to show that $y_1 < y_2$ implies
that $V(y_1) + V(y_2(y_1))$ is decreasing in $y_1$, or
equivalently,
\[
V_t(y_1) + V_t(y_2(y_1)) \frac{dy_2}{dy_1} < 0
\] 
where
$V(t)$ is the volume of the quoin, $V_t(t)$ is the
derivative of the volume, and $y_2$ is an implicit function
of $y_1$.

\begin{figure}
\includegraphics{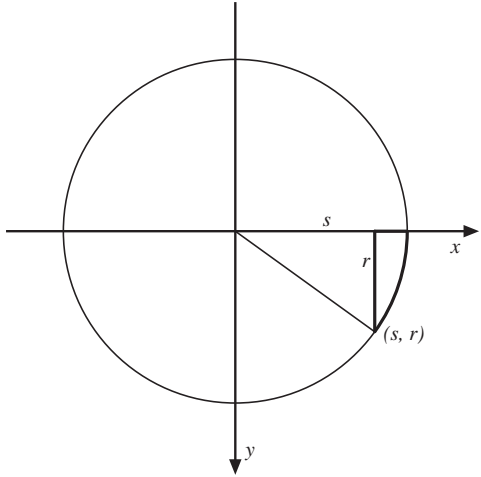}
%\centerline{ \psfig{file=quoin1.eps} }
\caption{Top view of quoin.}
\label{fig:quoin1}
\end{figure}

\begin{figure}
\includegraphics{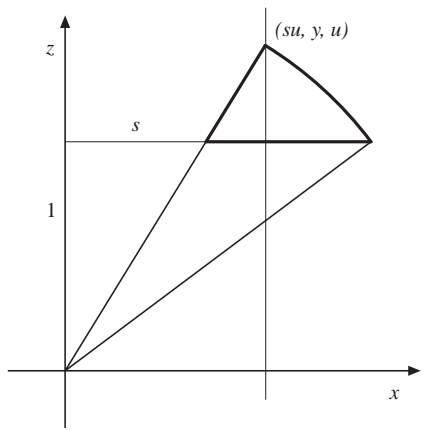}
%\centerline{ \psfig{file=quoin2.eps} }
\caption{Side view of quoin.}
\label{fig:quoin2}
\end{figure}

We construct the volume of a quoin by integrating the area of
a slice.  We place the quoin in a convenient coordinate system.
See Figures~\ref{fig:quoin1} and~\ref{fig:quoin2}.
The truncating sphere has equation $x^2 + y^2 + z^2 = 2$.
At the base of the quoin, $z=1$, so $x = \sqrt{1-y^2}$ gives
the location of the right-boundary of the quoin.
The plane forming the left face of the quoin is given by the
equation $x = s z$, so the ridge of the quoin is given by the
curve $(s u, y, u)$, where $u = \sqrt{\frac{2-y^2}{1+s^2}}$.

Hence the area of a slice parallel to the $x$-$z$ plane is given
by the formula
\[
A(t,y) = \frac{1}{2}(s u - s)(u - 1) + 
\int_{s u}^{\sqrt{1-y^2}}{(\sqrt{2-x^2-y^2}-1) \,dx}.
\]
The volume of a quoin is therefore given by the formula
\[
V(t) = \int_0^r{A(t,y)\,dy}.
\]

We actually only need to compute $V_t(t)$, 
which is fortunate,
since the explicit formula for $V(t)$ is somewhat complicated.
We have 
\[
V_t(t) = \int_0^r{A_t(t,y)\,dy} + A(t,r) r_t,
\]
but $A(t,r) = 0$, so 
\[
V_t(t) = \int_0^r{A_t(t,y)\,dy}.
\]
So in addition, we only need $A_t(t,y)$,
%\[
%A_t(t,y) = (\frac{s}{2}(u^2+1) - \sqrt{1-y^2} +
%\int_{s u}^{\sqrt{1-y^2}}{\sqrt{2-x^2-y^2}\,dx})_t,
%\]
\[
A_t(t,y) = (\frac{s}{2}(u^2+1) - \sqrt{1-y^2} +
\int_{\frac{t}{4}\sqrt{2-y^2}}^{\sqrt{1-y^2}}{\sqrt{2-x^2-y^2}\,dx})_t,
\]
so
\[
A_t(t,y) = (\frac{s}{2}(u^2+1))_t -
\sqrt{2-\frac{t^2}{16}(2-y^2) - y^2}\frac{1}{4}\sqrt{2-y^2}
\]
which simplifies to
\[
A_t(t,y) = \frac{8}{(16-t^2)^{3/2}} - \frac{2-y^2}{2\sqrt{16-t^2}}.
\]
Hence
\[
V_t(t) = \frac{8r}{(16-t^2)^{3/2}} - \frac{r}{\sqrt{16-t^2}} + 
\frac{r^3}{6\sqrt{16-t^2}},
\]
which simplifies to
\[
V_t(t) = \frac{-2\sqrt{2}(8-t^2)^{3/2}}{3(16-t^2)^2}.
\]

\subsubsection{The solid angle constraint}
Holding the solid angle fixed, $y_2$ is an implicit 
function of $y_1$.  We now
make that relation explicit.
Using formulas from \cite{sp1}, the solid angle constraint, 
\[\sol(2,2,2,y_1,y_2,d)=c,\]
where $c$
is a constant, becomes
\[
2 \arctan(\frac{\sqrt{\Delta}}{2a}) = c.
\]
Let $x_1 = y_1^2$, $x_2 = y_2^2$, and $b = d^2$.  Then
\[
\Delta = -4b^2 -4(x_1-x_2)^2 + b(x_1(8-x_2) + 8x_2),
\]
and
\[
a = 32 - d - x_1 - x_2.
\]
So
\[
\frac{-4b^2 -4(x_1-x_2)^2 + b(x_1(8-x_2) + 8x_2)}{(32 - d - x_1 - x_2)^2} = c_1.
\]
Therefore
\[
\frac{dx_2}{dx_1} = -\frac{(16-x_2)(x_2 + b - x_1)}{(16-x_1)(x_1 + b - x_2)},
\]
and
\[
\frac{dy_2}{dy_1} = \frac{y_1}{y_2} \frac{dx_2}{dx_1},
\]
hence
\[
\frac{dy_2}{dy_1} = -\frac{y_1(16-x_2)(x_2 + b - x_1)}{y_2(16-x_1)(x_1 + b - x_2)}.
\]

We return to the relation which we wish to prove, that $y_1 < y_2$ implies
\[
V_t(y_1) + V_t(y_2) \frac{dy_2}{dy_1} < 0.
\]
Note that all of the denominators are positive.
Therefore clearing the denominators, we
find that the desired relation
is equivalent to
\begin{eqnarray*}
-(8-x_1)^{3/2}(16-x_2)y_2(x_1 + b - x_2) + \\
(8-x_2)^{3/2}(16-x_1) y_1 (x_2 + b - x_1) & < & 0,
\end{eqnarray*}
or
\begin{eqnarray*}
(16-x_1)^2 x_1 (8-x_2)^3 (b - x_1 + x_2)^2 < \\
(16-x_2)^2 x_2 (8-x_1)^3 (b +x_1 - x_2)^2.
\end{eqnarray*}
If we define
\[
g(x_1,x_2) = (16-x_1)^2 x_1 (8-x_2)^3 (b - x_1 + x_2)^2,
\]
then the
desired inequality is equivalent to
$g(x_1,x_2) < g(x_2,x_1)$ for $x_1 < x_2$.
There are several ways to prove this monotonicity relation.
One is to prove that the
polynomial
\[
\frac{g(x_1,x_2)-g(x_2,x_1)}{8(x_1-x_2)}
\]
is positive for all allowable values for $x_1$, $x_2$, and $b$.  Unfortunately, the
resulting polynomial has degree $6$, so the verification is somewhat unwieldy,
although easy enough using interval methods.

A simpler method involves a factorization of $g$ into
$g_1$ and $g_2$.  We show that $g_1$ and $g_2$ each
satisfy the monotonicity relation, and the relation then follows for $g$.

Define 
\[
g_1(x_1,x_2) = (16-x_1) x_1 (8-x_2) (b-x_1+x_2),
\]
and
\[
g_2(x_1,x_2) = (16-x_1)(8-x_2)^2(b-x_1+x_2).
\]
Clearly
$g = g_1 g_2$.  We then construct the polynomials
\[
p_1 = \frac{g_1(x_1,x_2) - g_1(x_2,x_1)}{x_1-x_2}
\]
and
\[
p_2 = \frac{g_2(x_1,x_2) - g_2(x_2,x_1)}{x_1-x_2}.
\]

Simplifying $p_1$ and $p_2$, we find that
\begin{align*}
p_1 = & 128b - 128x_1 -8b x_1 + 8 x_1^2 - 128 x_2 \\
	& + 32 x_1 x_2 + b x_1 x_2 -x_1^2 x_2 + 8 x_2^2 - x_1 x_2^2
\end{align*}
and
\begin{align*}
p_2 = & -2048 + 192b + 320x_1 -16b x_1 -16x_1^2 +
	320x_2 \\ & - 16b x_2 -32 x_1 x_2 +b x_1 x_2 + 
	x_1^2 x_2 - 16x_2^2 + x_1 x_2^2.
\end{align*}
These polynomials are quadratic in $x_1$ and $x_2$, 
and linear in $b$.  The coefficient of $b$ in $p_1$
is
\[
128 -8x_1 -8x_2 + x_1 x_2.
\]
The coefficient of $b$ in $p_2$ is
\[
192 - 16x_1 -16x_2 + x_1 x_2.
\]
Both coefficients are positive for 
$x_1$ and $x_2$ in $[16/2.51^2,2.51^2]$.  Therefore, the
minimum values of $p_1$ and $p_2$ occur 
when $b$ is at a minimum, $b=8$.

The minimum value of each polynomial
for values of $x_1$ and $x_2$ in the range $[16/2.51^2,2.51^2]$
is now easily computed.
Making the appropriate computations, we find that each polynomial
is indeed positive.  Hence the desired relation follows.

\subsubsection{The simplification}
We now apply the reduction argument.  If
we are not careful about how we apply the argument, however, this
reduction could introduce some complications.

We begin with a squashed quad cluster with consecutive upper
edge lengths $(y_1,y_2,y_3,y_4)$ and diagonal $d$ adjacent to
the first two upper edges.

Recall that we chose the diagonal of the quad cluster to be the
shorter of the two possible diagonals.  We refer to the other
possible diagonal as the {\em cross-diagonal}.  
Recall that the reduction fixes the length of the diagonal.

If the length of the cross-diagonal does not drop to $2\sqrt{2}$
under the perturbation, we arrive at the configuration with
edge lengths $(y_1', y_1', y_2', y_2')$ with diagonal $d$.

If the length of the cross-diagonal does drop to $2\sqrt{2}$, then
stop the perturbation.  This gives a quadrilateral
$(y_1', y_2', y_3', y_4')$ with diagonal $2\sqrt{2}$.
Applying the perturbation to each half independently, we find
that the score of each half is maximized by the configuration
$(y_1'', y_1'', y_2'', y_2'')$ with diagonal $2\sqrt{2}$.
We verify the relation for this arrangement in
Calculation~\ref{obtuse:vor2}.

If the length of the cross-diagonal did not drop to $2\sqrt{2}$,
switch to the cross-diagonal and repeat the process.  If the
(new) cross-diagonal does not drop to $2\sqrt{2}$, we have
arrived at the configuration $(y,y,y,y)$ with diagonal $d'$.
Choose a new diagonal $d''$ to be the shorter of the two
possible diagonals.  We verify the desired relation for
this arrangement in Calculation~\ref{obtuse:vor}.

Finally, we make a few comments about extra constraints in
the verifications.

Since the score of a quad cluster
is non-positive, and $m (2\sol_0) - b \le 0$ where
$\sol_0 = \sol( 2,2,2,2,2,2\sqrt{2} )$, we need only consider
quad clusters for which the solid angle exceeds $2\sol_0$.

The maximum length of the diagonal is $2.51 \sqrt{2}$, since
otherwise the triangles in the quadrilateral would be obtuse,
forcing the cross-diagonal to be shorter than the diagonal.
This would contradict our
original choice of the shortest diagonal.

In Calculation~\ref{obtuse:vor}, we assume that $d$ is the shortest diagonal.
Adding this constraint directly is tedious, since the
formula for the cross-diagonal of the quad cluster is somewhat
complicated.  We apply a simpler but weaker constraint, 
that the diagonal $d$ of a planar quadrilateral
with edge lengths $(y,y,y,y)$ is shorter than $d'$, the other planar
diagonal.  The constraint $d \le d'$ gives the constraint $d^2 \le 2 y^2$.
Since the cross-diagonal of the quad cluster is shorter than the
cross-diagonal of the planar quadrilateral, this constraint is weaker.

\subsection{Mixed Quad Clusters}
By Proposition 4.1 of \cite{FKC}, the score of a mixed quad cluster is
less than $-1.04 \scorept$.  Therefore, we may discard all mixed
quad clusters.

\section{Numerical Considerations}
The verifications of the relations required in this paper 
appear intractable using
traditional methods.  Therefore, we use a relatively new proof technique,
interval arithmetic via floating-point computer calculations.

Most real numbers are not representable in computer floating-point
format.  However, floating-point intervals may be found which contain
any real number.  Although the magnitude of real numbers representable in
fixed-length floating-point format is finite, the format also provides
for $\pm \infty$, which allows for interval containment of all reals.  
These intervals may be
added, multiplied, etc., and the resulting intervals will contain
the result of the operation applied to the real numbers which they
represent.

Since floating-point arithmetic is not exact, interval arithmetic
conducted using floating-point arithmetic is not optimal, in the
sense that the interval resulting from an operation will usually
be larger than the true resultant interval, due to roundoff.
However, barring hardware or software errors (implementation errors,
not roundoff errors), floating-point
interval arithmetic, unlike floating-point arithmetic, is correct, in
the sense that it provides correct interval bounds on the value
of a computation, while floating-point arithmetic alone only 
provides an approximation to the correct value of a computation.
We may therefore use interval arithmetic to prove mathematical
results.  Floating-point arithmetic alone, in the absence of
rigorous error analysis, cannot constitute a proof.

We implement floating-point interval arithmetic routines via
the IEEE 754 Standard for floating-point arithmetic \cite{IEEE}.

Implementation of interval arithmetic is straightforward using
directed rounding.  In addition to arithmetic functions, we require
interval implementations of the
square root and arctangent functions.  Fortunately, the IEEE standard
provides the square root function.  However, the arctangent function
is somewhat problematic, since the
standard math libraries do not provide explicit error bounds for their
implementations of the arctangent function.  In theory, they
should provide an accuracy for the arctangent routine of
$0.7$ ulps, meaning that the error is less than one unit in the
last place.  I add interval padding of the form 
$[v-\epsilon, v+\epsilon]$,
where $v$ is the computed value, and $\epsilon = 2^{-49}$.
This should be sufficient to guarantee proper interval containment,
assuming that the library routines are correctly implemented.

Armed with standard interval arithmetic and interval arithmetic
implementations of {\em sqrt} and {\em arctan}, we can implement 
interval arithmetic versions of all
the special functions required for proving the sphere packing
relations.

Evaluating these functions on cells, we get bounds.  Unfortunately,
these bounds are not very good.  The bounds which we get from
interval versions of the partial derivative functions are even
worse.  This means that cells have to be very small before we
can draw conclusions about the signs of the partials.  These bad
bounds are due to the inherent nature of interval arithmetic--it
produces worst-case results by design.

These bad bounds increase the complexity of the verifications tremendously.
Some verifications, using these bounds, 
require the consideration of billions or trillions
of cells, or worse.  Therefore, we needed a method for producing
better bounds than those which direct interval methods could
provide.

The method which we eventually discovered is to use Taylor series.
We compute explicit second (mixed) partial bounds for the major
special functions, and use these bounds to produce very good
interval bounds.  These bounds are computed in 
Calculations \ref{secpar:dih:qr} through \ref{secpar:vor}.
%The beauty of the Taylor series method is that the error
%enters the calculation only once, instead of compounding with
%every arithmetic operation.
Essentially, the Taylor method postpones the error bound until the
end of the computation, eliminating the error bound explosion which
occurs with a straightforward interval method implementation.

%\appendix

\section{Calculations}
The following inequalities have been proved by computer using
interval methods.  Let $S = S(y) = S(y_1,\ldots,y_6)$ denote
a tetrahedron parametrized by the edge lengths $(y_1,\ldots,y_6)$.
In addition, we often parametrize by the squares of the edge
lengths, $(x_1,\ldots,x_6)$.

Recall from Section~\ref{sec:dabound} that $m = 0.3621$, $b = 0.49246$,
$\epsilon = 0.0739626$ and $b_c = 0.253095$.

%Moved Dimension Reduction to end, so numbering matches better.

\subsection{Quarters}

\begin{calc}
$\gma(S) \le 0$ for $y \in [2,2.51]^5[2.51,2\sqrt{2}]$.
\label{q:gma}
\end{calc}

\begin{calc}
$\octavor(S) \le 0$ for $y \in [2.51,2\sqrt{2}][2,2.51]^5$.
\label{q:octavor}
\end{calc}

Define the corner cell 
$C = [2,2 + 0.51/16]^5[2\sqrt{2} - (2\sqrt{2}-2.51)/16,2\sqrt{2}]$.
%[2, 2.031875]
%[2, 2.031875]
%[2, 2.031875]
%[2, 2.031875]
%[2, 2.031875]
%[2.80852542944955,  2.82842712474619]

\begin{calc}
$\frac{d}{dy_i}\vor(S) < 0$ for $i=1,\ldots,5$ and $y \in C$.
\label{q:partials}
\end{calc}

\begin{calc}
\[\vor(S) + 0.125(\eta(y_1,y_2,y_6)^2 - 2) \le 0\] 
for $y_3=y_4=y_5=2$ and $y \in C$, and
\[\vor(S) + 0.125(\eta(y_4,y_5,y_6)^2 - 2) \le 0\] 
for $y_1=y_2=y_3=2$ and $y \in C$.
\label{q:corner}
\end{calc}

\begin{calc}
$\vor(S) \le 0$ for $y \in [2,2.51]^5[2.51,2\sqrt{2}]$, $y \notin C$,
using dimension-reduction.
\label{q:off:dimred}
\end{calc}

\begin{calc}
\[\vor(S) \le 0\] for $y \in [2,2.51]^5[2.51,2\sqrt{2}]$, $y \notin C$
with \[\eta(y_1,y_2,y_6)^2=2 \text{~or~} \eta(y_4,y_5,y_6)^2=2,\] not using
dimension-reduction.
\label{q:off:bdry}
\end{calc}

\subsection{Quasi-regular Tetrahedra}

Define $C = [2,2.51]^6$, and recall
\[
a_1 = 0.3860658808124052, a_2 = 0.4198577862,
d_0 = 1.4674.\]

\begin{calc}
Either
\[\gma(S) \le a_1 \dih(S) - a_2\]
or
\[\gma(S) \le -0.52\scorept\]
for $y \in C$, using dimension-reduction.
\label{qr:dihrel}
\end{calc}

\begin{calc}
Either 
\[\gma(S) - a_1 \dih(S) \le 3.48 \scorept - 2 \pi a_1 + 4 a_2\]
or
\[\dih(S) < d_0\]
or
\[\gma(S) \le -0.52\scorept\]
for $y \in C$, using dimension-reduction.
\label{qr:dihcut}
\end{calc}

\begin{calc}
Either 
\[\gma(S) + m \sol(S) + \epsilon( \dih(S) - 
\frac{2\pi}{5}) - b_c \le 0\]
or
\[\dih(S) > d_0\]
or
\[\gma(S) \le -0.52\scorept\]
for $y \in C$, using dimension-reduction.
\label{qr:pentacap}
\end{calc}

\subsection{Flat Quad Clusters}

Define $I = [2,2.51]^5[2.51,2\sqrt{2}]$, and define the
corner cell
\[C = [2,2 + 0.51/16]^5[2\sqrt{2} - (2\sqrt{2}-2.51)/16,2\sqrt{2}].\]

\begin{calc}
Either 
\[\gma(S) + m \sol(S) \le b/2\]
or
\[\eta(y_1,y_2,y_6)^2 > 2\]
or
\[\eta(y_4, y_5, y_6)^2 > 2\]
or
\[\gma(S) \le -1.04\scorept\]
for $y \in I$, using dimension reduction.
\label{flat:gma:dimred}
\end{calc}

\begin{calc}
Either 
\[\gma(S) + m \sol(S) \le b/2\]
or
\[\eta(y_1,y_2,y_6)^2 = 2 \text{~with~} \eta(y_4, y_5, y_6)^2 \le 2,\]
or
\[\eta(y_4, y_5, y_6)^2 = 2 \text{~with~} \eta(y_1,y_2,y_6)^2 \le 2,\]
or
\[\gma(S) \le -1.04\scorept\]
for $y \in I$, not using dimension-reduction.
\label{flat:gma:bdry}
\end{calc}

\begin{calc}
$\frac{d}{dy_i}\vor(S) < 0$ for $i=1,2,3$ and $y \in C$.
\label{flat:partials}
\end{calc}

%This one is tricky.  Two cases:  either the side face, (1 2 6) is large,
%or the top face, (4 5 6) is large.  In the first case, 
%partials force eta(1 2 6)^2 = 2.  So solve for y1.

\begin{calc}
This computation is somewhat tricky, since the vor scoring
constraint depends on both faces.  The partial derivative
information allows us to assume $y_3=2$.  The rest of the
analysis depends on which face is assumed to be large.

If the $(y_1, y_2, y_6)$ face is large, the partial derivative
information allows us to assume that the face constraint
is tight, so $\eta(y_1,y_2,y_6)^2=2$.  Therefore we can solve
for $y_1$ in terms of $y_2$ and $y_6$.  We can apply partial
derivative information for $y_4$ and $y_5$.
In this case, we prove
\[\vor(S) + m \sol(S) \le b/2\]
for $y_3=2$, $y \in C$.

If the $(y_4,y_5,y_6)$ face is large, we may assume that
$y_1=y_2=2$.  We then prove
\[\vor(S) + m \sol(S) \le b/2\] 
or
\[\eta(y_4,y_5,y_6)^2 < 2\]
for $y_1=y_2=y_3=2$, $y \in C$.
\label{flat:corner}
\end{calc}

\begin{calc}
Either 
\[\vor(S) + m \sol(S) \le b/2,\]
or
\[\eta(y_1,y_2,y_6)^2 < 2 \text{~and~} \eta(y_4, y_5, y_6)^2 < 2,\]
or
\[\vor(S) \le -1.04\scorept\]
for $y \in I$, $y \notin C$, using dimension-reduction and partial
derivative information.
\label{flat:vor:dimred}
\end{calc}

\begin{calc}
Either 
\[\vor(S) + m \sol(S) \le b/2,\]
with
\[\eta(y_1,y_2,y_6)^2 = 2 \text{~or~} \eta(y_4, y_5, y_6)^2 = 2,\]
or
\[\vor(S) \le -1.04\scorept\]
for $y \in I$, $y \notin C$, not using dimension-reduction.
\label{flat:vor:bdry}
\end{calc}

\subsection{Octahedra}

\begin{calc}
%Since each quarter in an octahedron shares
%all edges except $y_4$ with an adjacent quarter,
%if one of these edges has length at least $2.2$, the
%adjacent quarter shares that property.  We demonstrate
%that if one of these edges is at least $2.2$, the
%score of the quarter does not exceed $-0.52 \scorept$,
%therefore the two quarters together may be discarded.
We prove $\score(S) \le -0.52 \scorept$,
for each (appropriately scored) upright quarter
with edge lengths in the cell
$[2.51,2\sqrt{2}][2.2,2.51][2,2.51]^4$.  
\label{octa:peel}
\end{calc}

\begin{calc}
Recall $c = 0.1533667634670977$, and
$d = 0.2265$.  We prove
\[\gma(S) + c \dih(S) \le d\]
or
\[\gma(S) \le -1.04 \scorept\]
for $y \in [2.51,2.716][2,2.2]^5$.  Note that for both
faces adjacent to the diagonal,
\[\max \eta^2 = \eta(2.2,2.2,2.716)^2 < 2,\]
so all quarters in this cell are compression-scored.
We make use of dimension-reduction.
\label{octa:cut}
\end{calc}

\begin{calc}
We prove 
\[\gma(S) + m \sol(S) + \alpha \dih(S)  \le \frac{b}{4} + 
	\alpha \frac{\pi}{2}\]
or
\[\gma(S) \le -1.04 \scorept\]
for all compression-scored
quarters $S(y)$, where $\alpha = 0.14$,
\[y \in [2.716,2\sqrt{2}][2,2.2]^2[2,2.51][2,2.2]^2.\]
We use dimension-reduction.
\label{octa:gma:dih}
\end{calc}

\begin{calc}
We prove 
\[\vor(S) + m \sol(S) + \alpha \dih(S)  \le \frac{b}{4} + 
	\alpha \frac{\pi}{2}\]
or
\[\vor(S) \le -1.04 \scorept\]
for all vor analytic-scored
quarters $S(y)$, where $\alpha = 0.14$,
\[y \in [2.716,2.81][2,2.2]^2[2,2.51][2,2.2]^2.\]
\label{octa:vor:dih}
\end{calc}

\begin{calc}
We prove 
\[\gma(S) + m \sol(S) + \alpha \dih(S) + 
	\beta x_1 \le \frac{b}{4} + 
	\alpha \frac{\pi}{2} + 8 \beta\]
or
\[\gma(S) \le -1.04 \scorept\]
for all compression-scored
quarters $S(y)$, where $\alpha = 0.054$, $\beta = 0.00455$,
$x_1=y_1^2$, and
\[y \in [2.81,2\sqrt{2}][2,2.2]^2[2,2.51][2,2.2]^2.\]  
We use some dimension-reduction.
\label{octa:gma:corr}
\end{calc}

\begin{calc}
We prove 
\[\vor(S) + m \sol(S) + \alpha \dih(S) + 
	\beta x_1 \le \frac{b}{4} + 
	\alpha \frac{\pi}{2} + 8 \beta\]
or
\[\vor(S) \le -1.04 \scorept\]
for all vor analytic-scored
quarters $S(y)$, where $\alpha = 0.054$, $\beta = -0.00455$,
$x_1=y_1^2$, and
\[y \in [2.81,2\sqrt{2}][2,2.2]^2[2,2.51][2,2.2]^2.\]
\label{octa:vor:corr}
\end{calc}

\subsection{Pure Voronoi Quad Clusters}

Recall $\sol_0$ denotes the solid angle of the tetrahedron
$(2,2,2,2,2,2\sqrt{2})$.

Define the corner cell $C = [2, 2 + 0.51/8]^5[2\sqrt{2}, 2.84]$.
We denote truncated Voronoi scoring by $\score$.
The constraint that the dividing face be acute translates
into $x_1 + x_2 - x_6 \ge 0$.  In each computation we apply
dimension-reduction.

We begin with the acute case.

\begin{calc}
We prove
\[\score(S) + m \sol(S) - b/2 \le 0\]
or
\[\sol(S) < \sol_0\]
or
\[x_1 + x_2 - x_6 < 0\]
or
\[\score(S) \le -1.04 \scorept\] for
$y \in [2,2.51]^5[2.84,4]$.
\label{acute:cut}
\end{calc}

\begin{calc}
We prove
\[\score(S) + m \sol(S) - b/2 \le 0\]
or
\[\sol(S) < \sol_0\]
or
\[x_1 + x_2 - x_6 < 0\]
or
\[\score(S) \le -1.04 \scorept\] 
for
$y \in [2,2.51]^5[2\sqrt{2},2.84]$ with $y \notin C$.
\label{acute:vor}
\end{calc}

\begin{calc}
We prove
\[\score(S) + m \sol(S) - b/2 \le 0\]
or
\[\sol(S) < \sol_0\]
or
\[x_1 + x_2 - x_6 < 0,\]
$y \in C$.
\label{acute:corner}
\end{calc}

Finally, we consider the obtuse case.

\begin{calc}
We prove
\[\score(S) + m \sol(S) - b/2 \le 0\]
or
\[\sol(S) < \sol_0\]
or
\[\score(S) \le -0.52 \scorept\] 
or
\[ 2 y^2 < d^2 \]
for a symmetric pure Voronoi quad cluster composed of two copies of $S$,
where \[S = (2,2,2,y,y,d),\] 
$y \in [4/2.51,2.51]$ and $d \in [2\sqrt{2}, 2.51 \sqrt{2}]$.
\label{obtuse:vor}
\end{calc}

\begin{calc}
We prove
\[\score(S_1) + \score(S_2) + m (\sol(S_1) + \sol(S_2)) - b \le 0\]
or
\[\score(S_1) + \score(S_2) \le -1.04 \scorept\] 
or
\[\sol(S_1) + \sol(S_2) < 2\sol_0\]
for a pure Voronoi quad cluster composed of two tetrahedrons $S_1$ and $S_2$,
where \[S_i = (2,2,2,y_i,y_i,2\sqrt{2}),\] 
$y_i \in [4/2.51,2.51]$.
\label{obtuse:vor2}
\end{calc}

\subsection{Dimension Reduction}

\begin{calc}
The polynomial derived for the dimension-reduction argument is
positive for $x \in [4,2.51^2]^6$ and $x \in [4,2.51^2]^5[4,8]$.
\label{d:dimred}
\end{calc}

\subsection{Second Partial Bounds}
We compute all second partials $\frac{d^2}{dx_i dx_j}$ in terms of
$x_i$, the squares of the edge lengths.  We do each computation
twice, once for \qrtets\ and once for quarters.  We compute the
second partials of $\dih$, $\sol$, $\gma$ volume, and $\vor$
volume (the vor analytic volume).  Since the scoring functions
are linear combinations of $\sol$ and the volume terms, we may
derive second partial bounds for $\gma$ and $\vor$ from these.

With the application of additional computer power, these
bounds could be improved.  These bounds were computed using
$16$ subdivisions.  While using $32$ subdivisions would improve
the bounds by a factor of $2$, perhaps, the time required
for the computations increases by a factor of $64$.

\begin{calc}
For \qrtets, the second partials of $\dih$ lie in
\[
[-0.0926959464,  0.0730008897].
%   1234567890     1234567890
\]
\label{secpar:dih:qr}
\end{calc}

\begin{calc}
For quarters, the second partials of $\dih$ lie in
\[
[-0.2384125007,   0.169150875].
%   1234567890      1234567890
\]
\label{secpar:dih}
\end{calc}

\begin{calc}
For \qrtets, the second partials of $\sol$ lie in
\[
[-0.0729140255,  0.088401996].
%   1234567890     1234567890
\]
\label{secpar:sol:qr}
\end{calc}

\begin{calc}
For quarters, the second partials of $\sol$ lie in
\[
[-0.1040074557,   0.1384785805].
%   1234567890      1234567890
\]
\label{secpar:sol}
\end{calc}

\begin{calc}
For \qrtets, the second partials of $\gma$ volume lie in
\[
[-0.0968945273,  0.0512553817].
%   1234567890     1234567890
\]
\label{secpar:gma:qr}
\end{calc}

\begin{calc}
For quarters, the second partials of $\gma$ volume lie in
\[
[-0.1362100221,   0.1016538923].
%   1234567890      1234567890
\]
\label{secpar:gma}
\end{calc}

\begin{calc}
For \qrtets, the second partials of $\vor$ volume lie in
\[
[-0.1856683356,   0.1350478467].
%   1234567890      1234567890
\]
\label{secpar:vor:qr}
\end{calc}

\begin{calc}
For quarters, the second partials of $\vor$ volume lie in
\[
[-0.2373892383,   0.1994181009].
%   1234567890      1234567890
\]
\label{secpar:vor}
\end{calc}

The computed $\gma$ second partials then lie in
\[
[-0.2119591984, 0.2828323141],
%   1234567890    1234567890
\]
for \qrtets\ and quarters.

Likewise, the computed $\vor$ second partials then lie in
\[
[-0.7137209962, 0.8691765157],
%   1234567890    1234567890
\]
for \qrtets\ and quarters.

\appendix

\section{Computer Code}
I would like to include the C code used in performing the verifications.
However, there is simply too much code to include here.  I will
make the code available electronically \cite{code}.

The code is not particularly beautiful.  This is due in part
to the fact that there is no ``natural'' representation for
interval arithmetic in C.  The ``operator overloading" 
available in C++ allows for prettier notation, but extracts
a serious performance penalty.  For example, there are much more
efficient methods for implementing interval arithmetic on
polynomials than the binary approach which operator 
overloading requires.  As the speed of execution is very
important, I chose more cumbersome notation instead of sleeker,
but significantly slower, code.  Incidentally, speed was the
reason why I chose not to use an external interval arithmetic
package.

There is a fair amount
of duplication in the code, since it evolved significantly
over time, as new methods eclipsed older ones.  For example,
some of the code is specialized to treat only tetrahedra
with acute faces, since this allows for some simplifications
which improve the efficiency of the routines.  With the
development of the Taylor method, these improvements became
less significant, but are still part of the code, mostly due
to inertia on my part.

Since each verification typically has some unique characteristics,
the code for conducting each verification has unique elements. 
This makes it difficult to represent everything compactly.

The code should be fairly portable.  
However, rounding control routines vary from
platform to platform, as do special commands, like FMADD, the
floating-point fused multiply-add command.  Therefore, some
minor changes will be required to localize the code to
a particular platform.

To give the flavor of the essentials of the code, I will
include code for computing the solid angle of a tetrahedron, 
and code for the
verification of Calculation~\ref{octa:gma:dih}.

\subsection{Computing the solid angle}
By Lemma~8.4.2 of \cite{sp1}, the formula for the solid angle
of a tetrahedron $S$ with edge lengths $(y_1,y_2,\ldots,y_6)$ is
given by
\[
\sol(S) = 2 \arctan(\frac{\sqrt{\Delta}}{2a}).
\]
Here
\begin{align*}
a(y_1,y_2,\ldots,y_6) = & y_1 y_2 y_3 + 
				\frac{1}{2} y_1 (y_2^2 + y_3^2 - y_4^2) \\
					& +	\frac{1}{2} y_2 (y_1^2 + y_3^2 - y_5^2) + 
					\frac{1}{2} y_3 (y_1^2 + y_2^2 -y_6^2)
\end{align*}
and
\begin{align*}
\Delta(x_1,\ldots,x_6) = & x_1 x_4 (-x_1 + x_2 + x_3 - x_4 + x_5 + x_6) \\
					& + x_2 x_5(x_1 - x_2 + x_3 + x_4 - x_5 + x_6)		\\
					& + x_3 x_6(x_1 + x_2 - x_3 + x_4 + x_5 - x_6)		\\
					& - x_2 x_3 x_4 - x_1 x_3 x_5 - x_1 x_2 x_6 - x_4 x_5 x_6,
\end{align*}
where $x_i = y_i^2$ for $i=1,\ldots,6$.

Note that the function $a$ is increasing in $y_1,y_2,$ and
$y_3$ on $[2,4]$ and is decreasing in the variables
$y_4, y_5,$ and $y_6$ on the same interval.  We use this
fact to simplify the interval calculation of $a$.

We represent an interval $t$ as $\big[\underline{t},\overline{t}\big]$.
In C, we represent this interval as a double-precision
floating-point array with two elements, {\tt double t[2]}.
We define $\underline{t} = \text{\tt t[0]}$ and
$\overline{t} = \text{\tt t[1]}$.

If each $y_i$ represents an interval, we denote the cell 
$(y_1,\ldots,y_6)$ in C as a double-precision floating-point
array with twelve elements, {\tt double y[12]}.  We define
$\underline{y_i} = \text{\tt y[2(i-1)]}$, and
$\overline{y_i} = \text{\tt y[2(i-1)+1]}$.

We begin with the implementation of $\Delta$.  The macros
{\tt ROUND\_DOWN} and {\tt ROUND\_UP} change the rounding
direction for floating-point computations.  We attempt
to take advantage of the polynomial structure, reducing
the interval overhead as much as possible.  

Essentially, we are taking advantage of the fact that it
is much more efficient to compute $a + b + c$ as
\[
[\underline{a} + \underline{b} + \underline{c}, 
	\overline{a} + \overline{b} + \overline{c}],
\]
since we only need to change the rounding direction once, 
compared to computing $(a + b) + c$, which would be required
by a binary interval arithmetic operator.  We therefore avoid
much of the function-call overhead, and in addition, we
won't cause as many floating-point pipeline stalls, since we
don't have to change the rounding direction as often.

{\footnotesize
\begin{verbatim}
void i_bigdelta( double x[12], double out[2] )
{
  double pterms[2], nterms[2], p1, p2, p3, p4;

  ROUND_DOWN;
  p1 = x[0]*x[6]*(x[2] + x[4] + x[8] + x[10]);
  p2 = x[2]*x[8]*(x[0] + x[4] + x[6] + x[10]);
  p3 = x[4]*x[10]*(x[0] + x[2] + x[6] + x[8]);
  pterms[0] = p1 + p2 + p3;
  p1 = x[0]*x[6]*(x[0] + x[6]);
  p2 = x[2]*x[8]*(x[2] + x[8]);
  p3 = x[4]*x[10]*(x[4] + x[10]);
  p4 = x[4]*(x[2]*x[6] + x[0]*x[8]) + x[10]*(x[0]*x[2] + x[6]*x[8]);
  nterms[0] = p1 + p2 + p3 + p4;
  ROUND_UP;
  p1 = x[1]*x[7]*(x[3] + x[5] + x[9] + x[11]);
  p2 = x[3]*x[9]*(x[1] + x[5] + x[7] + x[11]);
  p3 = x[5]*x[11]*(x[1] + x[3] + x[7] + x[9]);
  pterms[1] = p1 + p2 + p3;
  p1 = x[1]*x[7]*(x[1] + x[7]);
  p2 = x[3]*x[9]*(x[3] + x[9]);
  p3 = x[5]*x[11]*(x[5] + x[11]);
  p4 = x[5]*(x[3]*x[7] + x[1]*x[9]) + x[11]*(x[1]*x[3] + x[7]*x[9]);
  nterms[1] = p1 + p2 + p3 + p4;
  I_SUB( pterms, nterms, out );
  /* Only care about simplices for which delta > 0 */
  if( out[0] < 0.0 )
    out[0] = 0.0;
  if( out[1] < 0.0 )
    out[1] = 0.0;
}
\end{verbatim}
}
By Lemma~8.1.4 of \cite{sp1}, if $\Delta$ is negative, the simplex associated 
with the squares of the edge lengths
$(x_1,\ldots,x_6)$ does not exist.  
We therefore discard simplices with $\Delta < 0$.

We need to compute $\sqrt{\Delta}$.  In all of our code, we may
assume that the argument of the square root function is
non-negative.  If part of an interval is negative, we truncate
it so {\tt sqrt} won't complain.

{\footnotesize
\begin{verbatim}
void i_sqrt( double x[2], double out[2] )
{
  ROUND_DOWN;
  if( x[0] < 0.0 )
    out[0] = 0.0;
  else
    out[0] = sqrt( x[0] );
  ROUND_UP;
  out[1] = sqrt( x[1] );
}
\end{verbatim}
}

We next compute the minimum value of $a$.  We take advantage
of the monotonicity relation which we observed previously.
In the discussion following Lemma~8.4.2 of \cite{sp1}, we
find that $a>0$ for any simplex with edges in the
interval $[2,4]$.  This fact is
obvious for simplices for acute faces.  We may therefore
adjust the value of {\tt min\_a} should it return a
negative number.

{\footnotesize
\begin{verbatim}
/* min_a( y ) = a( y[0], y[2], y[4], y[7], y[9], y[11] ) */

double min_a( double y[12] )
{
  double y2[12];
  double p1, p2, p3, y123;
  int i;

  ROUND_DOWN;
  for( i=0; i<6; i+=2 )
    y2[i] = y[i]*y[i];
  ROUND_UP;
  for( i=7; i<12; i+=2 )
    y2[i] = y[i]*y[i];
  ROUND_DOWN;
  p1 = y2[2] + y2[4] - y2[7];
  p2 = y2[0] + y2[4] - y2[9];
  p3 = y2[0] + y2[2] - y2[11];
  y123 = y[0]*y[2]*y[4];
  p1 = y123 + 0.5*(y[0]*p1 + y[2]*p2 + y[4]*p3);
  if( p1 < 0.0 )
    p1 = 0.0;
  return( p1 );
}
\end{verbatim}
}

Often, we only need one half of the interval bound
for $a$, which is why I separated the routines.

{\footnotesize
\begin{verbatim}
/* max_a( y ) = a( y[1], y[3], y[5], y[6], y[8], y[10] ) */

double max_a( double y[12] )
{
  double y2[12];
  double p1, p2, p3, y123;
  int i;

  ROUND_DOWN;
  for( i=1; i<6; i+=2 )
    y2[i] = y[i]*y[i];
  ROUND_UP;
  for( i=6; i<12; i+=2 )
    y2[i] = y[i]*y[i];
  ROUND_UP;
  p1 = y2[3] + y2[5] - y2[6];
  p2 = y2[1] + y2[5] - y2[8];
  p3 = y2[1] + y2[3] - y2[10];
  y123 = y[1]*y[3]*y[5];
  p1 = y123 + 0.5*(y[1]*p1 + y[3]*p2 + y[5]*p3);
  return( p1 );
}
\end{verbatim}
}

We now compute the solid angle.  {\tt ATANERR} is a constant
which corrects for the potential error in computing
the arctangent.

{\footnotesize
\begin{verbatim}
void i_solid( double y[12], double sqrtdelta[2], 
  double out[2] )
{
  double max, temp, mina;
  
  temp = sqrtdelta[0];
  mina = 2.0*max_a( y );
  ROUND_DOWN;
  max = temp/mina;
  out[0] = 2.0*(atan( max ) - ATANERR);

  temp = sqrtdelta[1];
  mina = 2.0*min_a( y );
  ROUND_UP;
  max = temp/mina;
  out[1] = 2.0*(atan( max ) + ATANERR);
}
\end{verbatim}
}

\subsection{A sample verification}
The code for the verification of 
Calculation~\ref{octa:gma:dih} should
be representative of how the verifications
were conducted.

{
\footnotesize
\begin{verbatim}
/* octa_gma.c, by Samuel Ferguson, (c) 1997. */

#include "system_headers.h"
#include "sphere.h"
#include "interval.h"
#include "i_sphere.h"
#include "i_bounds.h"
#include "i_taylor.h"
#include "macros.h"


#define MAXDEPTH  32      /* Max recursion depth (was 10) */
#define PARTIALS  1       /* use partials, or not */

#define SUBDIV    2     /* (4)  number of subdivisions (1-d) */
#define SUBFRAC 0.5   /* (0.25)   1/SUBDIV (make this exact . . . */

#define MAXLEN    2.51
#define SCORECUT  -0.05758859149521   /* -1.04 pt */
#define CUTBD1    2.716
#define CUTBD2    2.81
#define Y4BOUND   TWO51_HI
#define PEELBD    2.2
#define ALPHA   0.14    /* 0.56/4 */ /* dihedral correction */

/* Even newer bounds:
line1[x_]:= 0.4922197796533495 - 0.3621*x
tetshift[x_]:= 0.253109 - 0.3621*x (eps = 0.07364)
*/

#define INTER1    0.4922197796533495
#define SLOPE1    0.3621

/* External variables */

extern double i_pi_const[2];
extern double i_pi_2_const[2];
extern double i_doct_const[2];
extern double i_two_pi_5_const[2];

/* External prototypes */

void i_init( void );

/* Global variables */

double cellcount;
double verifycount;
double partialcount;
double sccutcount;
double gmacount;
double vorcount;
double dualcount;
double celltotal;
int maxdepth;
int failed;
double slop;
double twosqrt2;

struct Cell_Info {
  int depth; 
  int cell_type;
  };

/* Prototypes */

void verifyoctagma( void );
void verify_cell( double tet[12], struct Cell_Info info );

/* Code */

void main( void )
{
  ROUND_NEAR;

  printf("Welcome to sphere packing routines, relation testing division.\n");
  verifyoctagma();
}


void verifyoctagma( void )
{
  int i, n;
  double tet[12];
  double diff, dt;
  time_t tstart, tstop;
  clock_t cstart, cstop;
  char *charptr;
  struct Cell_Info cell_info;
  
  printf("Welcome to OctaGma.  ");
#if PARTIALS
  printf("Using partials to reduce complexity.\n");
#else
  printf("Not using partials.\n");
#endif
  printf("Subdivisions: %d\n\n", SUBDIV);
  printf("Enter slop:  ");
  scanf("%lf", &slop);
  printf("\t\tSlop = %g\n", slop );
  
  i_init(); /* initialize interval constants */
  
  ROUND_DOWN;
  slop += 0.25*INTER1 + ALPHA*0.5*i_pi_const[0];
  ROUND_NEAR;

  time(&tstart);
  cstart = clock();
  
  charptr = ctime( &tstart );
  printf("Starting time:  \t");
  puts( charptr );
  
  celltotal = 0.0;
  for( n=1; n<6; n++ ) {
    printf("Starting case %d:\n", n);
    cellcount = 0.0;
    verifycount = 0.0;
    partialcount = 0.0;
    sccutcount = 0.0;
    vorcount = 0.0;
    gmacount = 0.0;
    dualcount = 0.0;
    maxdepth = 0;
    failed = 0;
    cell_info.depth = 0;
    cell_info.cell_type = -1;

    for( i=0; i<12; i++ )
      tet[i] = 2.0;
    tet[0] = CUTBD1;
    tet[1] = TWOSQRT2_HI;
    
    switch( n ) {
      case 1:
        tet[1] = CUTBD1;    /* tight */
        tet[7] = Y4BOUND;
        tet[9] = PEELBD;
        tet[11] = PEELBD;
        break;
      case 2:
        tet[1] = CUTBD1;    /* tight */
        tet[3] = PEELBD;
        tet[5] = PEELBD;
        tet[9] = PEELBD;
        tet[11] = PEELBD;
        break;
      case 3:
        tet[3] = PEELBD;
        tet[7] = Y4BOUND;
        tet[9] = PEELBD;
        break;
      case 4:
        tet[3] = PEELBD;
        tet[5] = PEELBD;
        tet[7] = Y4BOUND;
        break;
      case 5:
        tet[3] = PEELBD;
        tet[5] = PEELBD;
        tet[11] = Y4BOUND;
        break;

      default:
        printf("fell off end: this can't happen\n");
        break;
    }
    verify_cell( tet, cell_info );
    printf("cellcount    = %16.0f\t\t", cellcount);
    printf("maxdepth = %d\n", maxdepth);
    printf("sccutcount   = %16.0f\n", sccutcount);
    printf("gmacount     = %16.0f\t\t", gmacount);
    printf("vorcount     = %16.0f\n", vorcount);
    printf("dualcount    = %16.0f\t\t", dualcount);
    printf("verifycount  = %16.0f\n", verifycount);
    printf("partialcount = %16.0f\t\t", partialcount);
    diff = cellcount - verifycount - partialcount - sccutcount;
    printf("diff         = %16.0f\n", diff);
    if( !failed )
      printf("Verification succeeded.\n\n");
    else
      printf("MAXDEPTH exceeded.\n\n");
    celltotal += cellcount;
  } /* end case loop */
    
  ROUND_NEAR;
  time(&tstop);
  cstop = clock();
  printf("Done.\n");
  diff = ( (double) (cstop - cstart) )/CLOCKS_PER_SEC;
  charptr = ctime( &tstop );
  printf("Ending time:    \t");
  puts( charptr );
#if DA_SYSTEM == 2
  dt = ( double ) tstop - tstart;
#else
  dt = difftime( tstop, tstart );
#endif
  printf("Elapsed time: %g seconds\n", dt);
  printf("diff = %g\t\t", diff);
  printf("%g cells per second\n\n\n", celltotal/dt);
}


void verify_cell( double y[12], struct Cell_Info info )
{
  int i, j, bounds[6], n0, n1, n2, n3, n4, n5;
  int check, rel[5];
  double diff[6], t[12], val;
  double sc;
  double x[12], yp[12];
  double relconst[5], rel_val[2];
  double outvals[10];
  
  
  /* If part of the verification has already failed, bail. */
  if( failed )
    return;
  
  cellcount += 1.0;
  
  /* Compute square of edge lengths */
  
  ROUND_UP;
  for( i=1; i<12; i+=2 )
    x[i] = y[i]*y[i];
  ROUND_DOWN;
  for( i=0; i<12; i+=2 )
    x[i] = y[i]*y[i];
  if( x[1] > 8.0 )
    x[1] = 8.0;

  /* First try to verify over the current cell.  */
    
  /* Determine cell type */
  check = info.cell_type;
  if( check == -1 ) {  /* indeterminate */
    check = octa_scoring( x );
    info.cell_type = check;
  }
  
  if( check == 1 ) {
    gmacount +=1.0;
  }
  else if( check == 0 ) {
    vorcount += 1.0;
    return;
  }
  else
    dualcount += 1.0;
  
  /* Score cell properly */
  /* Order:  sol, gma, vor, octavor, dih */
  /*          0     1   2       3     4  */
  for( i=0; i<5; i++ ) {
    rel[i] = 0;
    relconst[i] = 0.0;
  }
  rel[1] = 1;
  relconst[0] = SLOPE1;
  relconst[1] = 1.0;
  relconst[4] = ALPHA;
  
  /* val = sc + SLOPE1*sph + ALPHA*dih */
  t_composite( rel, relconst, y, x, rel_val, t, 
    outvals );
  
  val = rel_val[1];
  sc = outvals[3];
  
  
  /* Discard cells with low score */
  if( sc < SCORECUT ) {
    sccutcount++;
    return;
  }
  
  /* Attempt verification */  
  if( val < slop ) {
    verifycount += 1.0;
    return;
  }
  
  /* if verification fails, try to bump off, or at least
    reduce dimension */

  for( i=0; i<12; i++ )
    yp[i] = y[i]; /* copy y */
  
#if PARTIALS
  
  /* Score cell properly (find min score) */
  if( check == 1 && outvals[2] > SCORECUT ) {
    /* t contains the relation partials */
    ROUND_NEAR;
    j = 0;
    for( i=0; i<6; i++ ) {
      diff[i] = yp[j+1] - yp[j];  /* compute differentials */
      j += 2;
    }
    /* Long edge is the first one, don't do that here. */
    /* Do edges 2 and 3 */
    j = 1;
    for( i=2; i<6; i+=2 ) {
      n0 = i + 1;
      if( diff[j] > 0.0 ) { /* ignore tight spots */
        if( t[n0] < 0.0 ) {
          if( yp[i] != 2.0 ) {
            partialcount += 1.0;
            return; /* bumped off */
          }
          else
            yp[n0] = 2.0; /* reduced dimension */
        }
        else {
          if( t[i] > 0.0 ) {
            if( yp[n0] != PEELBD ) {
              partialcount += 1.0;
              return; /* bumped off */
            }
            else
              yp[i] = PEELBD; /* reduced dimension */       
          }
        }
      }
      j++;
    }
    /* Do edges 5 and 6 */
    j = 4;
    for( i=8; i<12; i+=2 ) {
      n0 = i + 1;
      if( diff[j] > 0.0 ) { /* ignore tight spots */
        if( t[n0] < 0.0 ) {
          if( yp[i] != 2.0 ) {
            partialcount += 1.0;
            return; /* bumped off */
          }
          else
            yp[n0] = 2.0; /* reduced dimension */
        }
        else {
          if( t[i] > 0.0 ) {
            if( yp[n0] != PEELBD ) {
              partialcount += 1.0;
              return; /* bumped off */
            }
            else
              yp[i] = PEELBD; /* reduced dimension */       
          }
        }
      }
      j++;
    }
    /* Do edge 4 */
    j = 3;
    i = 6;
    n0 = 7;
    if( diff[j] > 0.0 ) { /* ignore tight spots */
      if( t[n0] < 0.0 ) {
        if( yp[i] != 2.0 ) {
          partialcount += 1.0;
          return; /* bumped off */
        }
        else
          yp[n0] = 2.0; /* reduced dimension */
      }
      else {
        if( t[i] > 0.0 ) {
          if( yp[n0] != Y4BOUND ) {
            partialcount += 1.0;
            return; /* bumped off */
          }
          else
            yp[i] = Y4BOUND;  /* reduced dimension */       
        }
      }
    }
    /* Now consider the long edge.  */
    /*
    j = 0;
    i = 0;
    n0 = 1;
    */
    if( diff[0] > 0.0 ) { /* ignore tight spots */
      if( t[1] < 0.0 ) {
        if( yp[0] != CUTBD1 ) {
          partialcount += 1.0;
          return; /* bumped off */
        }
        else
          yp[1] = CUTBD1; /* reduced dimension */
      }
      else {
        if( t[0] > 0.0 ) {
          if( yp[1] != TWOSQRT2_HI ) {
            partialcount += 1.0;
            return; /* bumped off */
          }
          else
            yp[0] = TWOSQRT2_HI; /* reduced dimension */        
        }
      }
    }
  }
#endif
  if( info.depth < MAXDEPTH ) {
    /* If all else fails, subdivide.  */
    info.depth++;
    if( info.depth > maxdepth )
      maxdepth = info.depth;
    ROUND_NEAR;
    j = 0;
    for( i=0; i<6; i++ ) {  /* scaled differentials */
      diff[i] = SUBFRAC*(yp[j+1] - yp[j]);
      j += 2;
      if( diff[i] > 0.0 )
        bounds[i] = SUBDIV;
      else
        bounds[i] = 1;
    }
    for( n0=0; n0<bounds[0]; n0++ ) {
      val = yp[0];
      t[0] = val + n0*diff[0];
      t[1] = val + (n0+1)*diff[0];
      for( n1=0; n1<bounds[1]; n1++ ) {
        val = yp[2];
        t[2] = val + n1*diff[1];
        t[3] = val + (n1+1)*diff[1];
        for( n2=0; n2<bounds[2]; n2++ ) {
          val = yp[4];
          t[4] = val + n2*diff[2];
          t[5] = val + (n2+1)*diff[2];
          for( n3=0; n3<bounds[3]; n3++ ) {
            val = yp[6];
            t[6] = val + n3*diff[3];
            t[7] = val + (n3+1)*diff[3];
            for( n4=0; n4<bounds[4]; n4++ ) {
              val = yp[8];
              t[8] = val + n4*diff[4];
              t[9] = val + (n4+1)*diff[4];
              for( n5=0; n5<bounds[5]; n5++ ) {
                val = yp[10];
                t[10] = val + n5*diff[5];
                t[11] = val + (n5+1)*diff[5];
                verify_cell( t, info );
                ROUND_NEAR;
              }
            }
          }
        }
      }
    }
  } else {  /* exceeded MAXDEPTH */
    failed = 1;
    ROUND_NEAR;
    printf("Cell failed:\n");
    printf("val = %g\n", val - slop);
    printf("sc = [%g, %g]\n", 
      outvals[2], outvals[3]);
    for( i=0; i<12; i+=2 )
      printf("%0.18f\t%0.18f\n", y[i], y[i+1]);
  }
}

\end{verbatim}
}

\section{Results}
The output of each verification contains information
about the number of cells considered, how the cells were
scored, how many cells were discarded and why, and the
number of recursive subdivisions.  I will
make the results of each computation available electronically
\cite{code}.  As a representative sample, I include
the output from Calculation~\ref{octa:gma:dih}.

This calculation has several cases, which arise from
the dimension-reduction argument.  It is a
relatively simple calculation, since only several
thousand cells had to be considered.  The inequality
which we prove is tighter than that exhibited (and
required) earlier in the paper.  
The ``slop'' value refers to the amount of
relaxation from the ideal bound.

{
\footnotesize
\begin{verbatim}
Welcome to sphere packing routines, relation testing division.
Welcome to OctaGma.  Using partials to reduce complexity.
Subdivisions:   2

Enter slop:  1.0e-5
        Slop = 1e-05
Starting time:      Sun Jun 29 11:03:12 1997

Starting case 1:
cellcount    =              363     maxdepth = 5
sccutcount   =              137
gmacount     =              363     vorcount     =                0
dualcount    =                0     verifycount  =              180
partialcount =                0     diff         =               46
Verification succeeded.

Starting case 2:
cellcount    =              609     maxdepth = 4
sccutcount   =               66
gmacount     =              609     vorcount     =                0
dualcount    =                0     verifycount  =              504
partialcount =                0     diff         =               39
Verification succeeded.

Starting case 3:
cellcount    =             3959     maxdepth = 9
sccutcount   =             1398
gmacount     =             3406     vorcount     =                0
dualcount    =              553     verifycount  =             2297
partialcount =                4     diff         =              260
Verification succeeded.

Starting case 4:
cellcount    =             3203     maxdepth = 9
sccutcount   =             1150
gmacount     =             2738     vorcount     =                0
dualcount    =              465     verifycount  =             1838
partialcount =                1     diff         =              214
Verification succeeded.

Starting case 5:
cellcount    =             2025     maxdepth = 9
sccutcount   =               39
gmacount     =             1366     vorcount     =              144
dualcount    =              515     verifycount  =             1701
partialcount =                1     diff         =              284
Verification succeeded.

Done.
Ending time:        Sun Jun 29 11:03:20 1997

Elapsed time: 8 seconds
diff = 7.83333      1269.88 cells per second


\end{verbatim}
}

Some verifications require the consideration of millions
of cells.  Without the use of dimension-reduction and
Taylor methods, they could
require billions of cells, or could even lie beyond our
current computer resources.

\end{document}